\documentclass{article}
\usepackage{graphicx}
\usepackage{amsmath}
\usepackage{amsfonts}
\usepackage{amssymb}
\newtheorem{theorem}{Theorem}

\newtheorem{corollary}[theorem]{Corollary}

\newtheorem{definition}[theorem]{Definition}
\newtheorem{example}[theorem]{Example}

\newtheorem{lemma}[theorem]{Lemma}

\newtheorem{proposition}[theorem]{Proposition}
\newtheorem{property}[theorem]{Property}
\newtheorem{remark}[theorem]{Remark}

\newenvironment{proof}[1][Proof]{\textbf{#1.} }{\ \rule{0.5em}{0.5em}}

\begin{document}

\title{Reconsideration of the multivariate moment problem and  a new method for
approximating multivariate integrals}
\author{Ognyan Kounchev and Hermann Render}
\maketitle
\begin{abstract}
Due to its intimate relation to Spectral Theory and Schr\"{o}dinger operators,
the multivariate moment problem has been a subject of many researches, so far
without essential success (if one tries to compare with the one--dimensional
case). In the present paper we reconsider a basic axiom of the standard
approach - the positivity of the measure. We introduce the so--called
pseudopositive measures instead. One of our main achievements is the solution
of the moment problem in the class of the pseudopositive measures. A measure
\ $\mu$ is called pseudopositive if its Laplace-Fourier coefficients
$\mu_{k,l}\left(  r\right)  ,$ $r\geq0,$ in the expansion in spherical
harmonics are non--negative. Another main profit of our approach is that for
pseudopositive measures we may develop efficient ''cubature formulas'' by
generalizing the classical procedure of Gauss--Jacobi: for every integer
\ $p\geq1$ we construct a new pseudopositive measure \ $\nu_{p}$ having
''minimal support'' and such that $\mu\left(  h\right)  =\nu_{p}\left(
h\right)  $ for every polynomial $h$ with \ $\Delta^{2p}h=0.$ The proof of
this result requires application of the famous theory of Chebyshev, Markov,
Stieltjes, Krein for extremal properties of the Gauss-Jacobi measure, by
employing the classical orthogonal polynomials $p_{k,l;j},$ $j\geq0,$  with
respect to every measure $\mu_{k,l}.$ As a byproduct we obtain a notion of
multivariate orthogonality defined by the polynomials $p_{k,l;j}$. A major
motivation for our investigation has been the further development of new
models for the multivariate Schr\"{o}dinger operators, which generalize the
classical result of M. Stone saying that the one--dimensional orthogonal
polynomials represent a model for the self--adjoint operators with simple
spectrum. 
\end{abstract}

\section{Introduction}

The univariate moment problem is one of the cornerstones of Mathematical
Analysis where several areas of Pure and Applied Mathematics meet -- continued
fractions, quadrature formulas, orthogonal polynomials, analytic functions,
finite differences, operator and spectral theory, scattering theory and
inverse problems, probability theory, and last but not least, control theory,
see e.g. the collection of surveys in \cite{landau} and the comprehensive
recent account \cite{simon} on the numerous applications of the moment problem
to spectral theory. On the other hand, the multivariate case is much more
complicated, and we refer to \cite{BCR84}, \cite{CuFi00}, \cite{McGr80},
\cite{PuVa99}, \cite{Schm91a}, \cite{StSz98} and the references given there
for some recent developments. However, the state of the art in the
multivariate moment problem seems to be well characterized by a remark in the
versatile survey \cite[p. 47]{Fugl83}, saying that only comparatively little
from the comprehensive theory of the classical moment problem has been
extended to dimension $d>1.$

The main purpose of the present paper is to introduce a \emph{modified}
\emph{moment problem} for which the solutions are in general signed measures
and belong to the class of what we call \emph{pseudo--positive measures}. The
motivation for this new notion is the possibility to generalize the univariate
Gau\ss--Jacobi quadratures to the multivariate setting, and thus to
approximate multivariate integrals in a new stable way. Let us emphasize that
we do not claim to solve the multivariate moment problem in its classical formulation.

In order to make our approach clear, let us first recall the usual formulation
of the multivariate moment problem: it asks for conditions on a sequence of
real numbers $c=\left\{  c_{\alpha}\right\}  _{\alpha\in\mathbb{N}_{0}^{d}}$
(here $\mathbb{N}_{0}$ denotes the set of all non--negative integers and we
use the multi--index notation $x^{\alpha}=x_{1}^{\alpha_{1}}x_{2}^{\alpha_{2}%
}\cdot\cdot\cdot x_{d}^{\alpha_{d}}$), such that there exists a
\emph{non-negative} measure $\mu$ on $\mathbb{R}^{d}$ with
\begin{equation}
c_{\alpha}=\int_{\mathbb{R}^{d}}x^{\alpha}d\mu\left(  x\right)
\label{eqmoment}%
\end{equation}
for all $\alpha\in\mathbb{N}_{0}^{d}.$ Let us denote by $\mathbb{C}\left[
x_{1},x_{2},...,x_{d}\right]  $ the space of all polynomials in $d$ variables
with complex coefficients. With the sequence $c=\left\{  c_{\alpha}\right\}
_{\alpha\in\mathbb{N}_{0}^{d}}$ we associate by linear extension to
$\mathbb{C}\left[  x_{1},x_{2},...,x_{d}\right]  $ a\textbf{\ }functional
$T_{c}$, by putting
\begin{equation}
T_{c}\left(  x^{\alpha}\right)  :=c_{\alpha}\qquad\text{for }\alpha
\in\mathbb{N}_{0}^{d}. \label{eqTc}%
\end{equation}
By a theorem of Haviland, a necessary and sufficient condition for the
existence of a non-negative measure $\mu$ satisfying (\ref{eqmoment}) is the
\emph{positivity of the sequence} $c=\left\{  c_{\alpha}\right\}  _{\alpha
\in\mathbb{N}_{0}^{d}}$, i.e. $P\geq0$ implies $T_{c}\left(  P\right)  \geq0$
for all $P\in\mathbb{C}\left[  x_{1},x_{2},...,x_{d}\right]  $ (here $P\geq0$
means that $P\left(  x\right)  \geq0$ for all $x\in\mathbb{R}^{d}$ ), cf.
\cite[p. 111]{Berg1987}. As is well known, the condition of positivity is
difficult to apply, and in practice one uses the weaker and easier to check
algebraic condition that the sequence $c=\left\{  c_{\alpha}\right\}
_{\alpha\in\mathbb{N}_{0}^{d}}$ is \emph{positive definite; }by definition the
sequence $\left\{  c_{\alpha}\right\}  _{\alpha\in\mathbb{N}_{0}^{d}}$ is
\emph{positive definite}\footnote{Some authors \cite{Akhi65} call such
sequences ''positive'', while we are closer to the terminology of say
\cite{BCR84}.} if and only if
\[
T_{c}\left(  P^{\ast}P\right)  \geq0\qquad\text{ for all }P\in\mathbb{C}%
\left[  x_{1},x_{2},...,x_{d}\right]  ;
\]
here $P^{\ast}$ is the polynomial whose coefficients are the complex
conjugates of the coefficients of $P.$ Let us remind that for $d=1,$
positivity and\ positive--definiteness of $T_{c}$\ are equivalent. However,
for $d>1$ this is \emph{not true}, and this is a consequence of the fact,
already known to D. Hilbert, that there exist non--negative polynomials which
are not sums of squares of other polynomials, see e.g. \cite{Gelfand},
\cite{BCR84}.

As mentioned above, we shall consider a different setting of the moment
problem. Let us postpone at the moment the motivation for our approach, and
let us concentrate on our new setting which needs some technical preparations
from the theory of harmonic functions. We assume that for each $k=0,1,2,...,$
the functions $Y_{k,l}:\mathbb{R}^{d}\rightarrow\mathbb{R},$ $l=1,...,a_{k},$
form a basis of the set of all harmonic homogeneous complex-valued
polynomials\footnote{For a reader not familiar with the spherical harmonics it
will be enough to consider the two--dimensional case $d=2$ where the basis is
simple and given in Section \ref{S2d}. On the other hand in Section
\ref{Sstrip} we develop our theory in the case of the strip where all what one
needs is expansion in Fourier series.} of degree $k\in\mathbb{N}_{0},$ and
they are orthonormal with respect to the scalar product $\left\langle
f,g\right\rangle _{\mathbb{S}^{d-1}}:=\int_{\mathbb{S}^{d-1}}f\left(
\theta\right)  \overline{g\left(  \theta\right)  }d\theta,$ where
$\mathbb{S}^{d-1}:=\left\{  x\in\mathbb{R}^{d}:\left|  x\right|  =1\right\}  $
is the unit sphere, and $r=\left|  x\right|  =\sqrt{x_{1}^{2}+....+x_{d}^{2}}$
is the euclidean norm, cf. \cite{ABR92} or \cite{StWe71}. For $x\in
\mathbb{R}^{d}$ we will use further the representation
\[
x=r\theta\qquad\text{for }r\geq0,\ \theta\in\mathbb{S}^{d-1}.
\]
The functions $Y_{k,l}$ are called \emph{solid harmonics} and their
restrictions to $\mathbb{S}^{d-1}$ \emph{spherical harmonics}. An important
property of the system $\left|  x\right|  ^{2j}Y_{k,l}\left(  x\right)  ,$
$j,k\in\mathbb{N}_{0},$ $l=1,...,a_{k}$, is that it forms a\textbf{\ }%
\emph{basis for} $\mathbb{C}\left[  x_{1},x_{2},...,x_{d}\right]  $. This
result follows from the Gau\ss\ decomposition of polynomials which says (cf.
\cite[Theorem 5.6, Theorem 5.21, p. 77 and p. 90]{ABR92}, \cite{Sob}, or
\cite[Theorem 10.2]{Koun00}) that every polynomial $P$ may be expanded in the
following way,
\begin{equation}
P\left(  x\right)  =\sum_{j=0}^{\left[  \deg\left(  P\right)  /2\right]
+1}\left|  x\right|  ^{2j}p_{j}\left(  x\right)  , \label{Gaussdecomposition}%
\end{equation}
where $p_{j}$ are harmonic polynomials, $\deg P$ denotes the degree of $P$ and
$\left[  x\right]  $ is the integer part of a real number $x.$ Since each
$p_{j}$ is a linear combination of the solid harmonics $Y_{k,l}\left(
x\right)  ,k\in\mathbb{N}_{0},\ l=1,2,...,a_{k},$ it is clear that the system
$\left|  x\right|  ^{2j}Y_{k,l}\left(  x\right)  ,$ $j,k\in\mathbb{N}_{0},$
$l=1,...,a_{k}$, is a basis for $\mathbb{C}\left[  x_{1},x_{2},...,x_{d}%
\right]  $. It is instructive to discuss the relationship between the
Gau\ss\ decomposition and the Laplace-Fourier series: recall that for a
sufficiently nice function $f:\mathbb{R}^{n}\rightarrow\mathbb{C}$ (e.g.
continuous) the expansion
\begin{equation}
f\left(  r\theta\right)  =\sum_{k=0}^{\infty}\sum_{l=1}^{a_{k}}f_{k,l}\left(
r\right)  Y_{k,l}\left(  \theta\right)  , \label{frtheta}%
\end{equation}
is the \emph{Laplace--Fourier series }with the \emph{Laplace--Fourier
coefficients} given by
\begin{equation}
f_{k,l}\left(  r\right)  =\int_{\mathbb{S}^{d-1}}f\left(  r\theta\right)
Y_{k,l}\left(  \theta\right)  d\theta. \label{eqfourier}%
\end{equation}
Suppose now that $f$ is a polynomial: then (\ref{Gaussdecomposition}) implies
that the Laplace-Fourier series (\ref{frtheta}) is a finite series, and the
functions $f_{k,l}\left(  r\right)  r^{-k}$ are \emph{polynomials} in the
variable $r^{2}.$ Moreover from (\ref{Gaussdecomposition}) directly follows
that each $f_{k,l}\left(  r\right)  r^{-k},k\in\mathbb{N}_{0}%
,\ l=1,2,...,a_{k},$ is a polynomial of degree $\leq2s-2$ if and only if for
all $x\in\mathbb{R}^{n}$
\begin{equation}
\Delta^{s}f\left(  x\right)  =0. \label{delta0}%
\end{equation}
Here $\Delta$ denotes the Laplace operator defined by $\Delta=\frac
{\partial^{2}}{\partial x_{1}^{2}}+...+\frac{\partial^{2}}{\partial x_{d}^{2}%
},$ and $\Delta^{s}$ is the $s$-th iterate of $\Delta$ for integers $s\geq1.$
In view of (\ref{delta0}) let us recall that a function $f$ defined on an open
subset $U$ in $\mathbb{R}^{d}$ is \emph{polyharmonic of order} $s$ if
$\Delta^{s}f\left(  x\right)  =0$ for all $x\in U,$ see \cite{ACL83}%
.\footnote{In the last section we provide some remarks on the significance of
the polyharmonic functions in approximation theory which have motivated also
the present research.}

Now we come to the cornerstone of our approach. Let $T_{c}:$ $\mathbb{C}%
\left[  x_{1},x_{2},...,x_{d}\right]  \rightarrow\mathbb{C}$ be a functional
associated to a sequence of moments $c_{\alpha},\alpha\in\mathbb{N}_{0}^{d}.$
Using the basis $\left|  x\right|  ^{2j}Y_{k,l}\left(  x\right)
:j,k\in\mathbb{N}_{0},\ l=1,...,a_{k}$ one can define a problem equivalent to
the usual one (\ref{eqmoment}), by means of the sequence $\{c_{j}^{\left(
k,l\right)  }\}_{j\in\mathbb{N}_{0}}$ defined as follows,
\begin{equation}
c_{j}^{\left(  k,l\right)  }:=T_{c}\left(  \left|  x\right|  ^{2j}%
Y_{k,l}\left(  x\right)  \right)  \qquad\text{for }j=0,1,2,...
\label{defpseudoc}%
\end{equation}
In order to distinguish them from the usual moments $c_{\alpha},$ the numbers
$c_{j}^{\left(  k,l\right)  }$ are sometimes called \emph{distributed moments,
}see \cite{butkovskii}, \cite{butkovskii2}, \cite{kounchev84},
\cite{kounchev85}, \cite{kounchev87}, \cite{kounchev93}.

We say that the sequence $\left\{  c_{\alpha}\right\}  _{\alpha\in
\mathbb{N}_{0}^{d}}$ or the associated functional $T_{c}$ is
\emph{pseudo-positive definite} if for every fixed pair of indices $\left(
k,l\right)  $ with $k\in\mathbb{N}_{0}$ and $l=1,...,a_{k}$ the sequences
$\{c_{j}^{\left(  k,l\right)  }\}_{j\in\mathbb{N}_{0}}$ and $\{c_{j+1}%
^{\left(  k,l\right)  }\}_{j\in\mathbb{N}_{0}}$ are \emph{positive definite}.
Equivalently, for each solid harmonic $Y_{k,l}\left(  x\right)  ,k\in
\mathbb{N}_{0},\ l=1,...,a_{k},$ the \emph{component functional }%
$T_{k,l}:\mathbb{C}\left[  x_{1}\right]  \rightarrow\mathbb{C}$ defined by
\begin{equation}
T_{k,l}\left(  p\right)  :=T_{c}\left(  p(\left|  x\right|  ^{2}%
)Y_{k,l}\left(  x\right)  \right)  \text{ for }p\in\mathbb{C}\left[
x_{1}\right]  \label{defpseudo}%
\end{equation}
has the property that $T_{k,l}\left(  p^{\ast}\left(  t\right)  p\left(
t\right)  \right)  \geq0$ and $T_{k,l}\left(  tp^{\ast}\left(  t\right)
p\left(  t\right)  \right)  \geq0$ for all $p\in\mathbb{C}\left[
x_{1}\right]  .$

Now we can formulate and successfully solve the following modified moment
problem: \emph{Given a pseudo-positive definite sequence }$c=\left\{
c_{a}\right\}  _{\alpha\in\mathbb{N}_{0}^{d}}$\emph{\ and its associated
functional }$T_{c},$\emph{\ find the conditions for the existence of a signed
measure }$\mu$\emph{\ on }$\mathbb{R}^{d}$\emph{\ such that}
\begin{equation}
\int_{\mathbb{R}^{n}}P\left(  x\right)  d\mu=T_{c}\left(  P\right)  \text{ for
all }P\in\mathbb{C}\left[  x_{1},x_{2},...,x_{d}\right]  . \label{eqpoll}%
\end{equation}
Two remarks are important: first, we allow $\mu$ to be a \emph{signed} measure
on $\mathbb{R}^{n}$, and this requirement is motivated by our constructive
formulas for approximating integrals developed in later sections. Secondly, it
follows from (\ref{eqpoll}) that the measure $\mu,$ considered as a functional
on $\mathbb{C}\left[  x_{1},x_{2},...,x_{d}\right]  ,$ is pseudo-positive
definite. The remarkable thing which will be seen from our further development
is that problem (\ref{eqpoll}) has a solution $\mu$ which
is\emph{\ pseudo-positive} which means that the inequality
\begin{equation}
\int_{\mathbb{R}^{d}}h\left(  \left|  x\right|  \right)  Y_{k,l}\left(
x\right)  d\mu\left(  x\right)  \geq0 \label{defpspos}%
\end{equation}
holds for every non-negative continuous function $h:\left[  0,\infty\right)
\rightarrow\left[  0,\infty\right)  $ with compact support and for all pairs
of indices $\left(  k,l\right)  $ with $k\in\mathbb{N}_{0}$ and
$l=1,2,...,a_{k}.$

As a first evidence that the pseudo--positivity is a reasonable generalization
of the univariate positivity notion, we present in Section
\ref{Smomentproblem} the following \emph{solution to the modified moment
problem }(\ref{eqpoll}), provided in two steps: 1. By a classical
one--dimensional argument, for the component functionals $T_{k,l}$ associated
with the pseudo--positive definite functional $T_{c}:\mathbb{C}\left[
x_{1},x_{2},...,x_{d}\right]  \rightarrow\mathbb{C}$ there exist non--negative
univariate representing measures $\mu_{k,l}$ on $\left[  0,\infty\right)  .$
2. If they satisfy the \emph{summability assumption}
\begin{equation}
\sum_{k=0}^{\infty}\sum_{l=1}^{a_{k}}\int_{0}^{\infty}r^{N}r^{-k}d\mu
_{k,l}\left(  r\right)  <\infty\text{ for all }N\in\mathbb{N}_{0}
\label{maincond}%
\end{equation}
then there exists a pseudo--positive signed measure $\mu$ on $\mathbb{R}^{d}$
representing $T_{c}$, i.e. (\ref{eqpoll}) holds. Further the following
important identity
\begin{equation}
\int_{\mathbb{R}^{n}}f\left(  x\right)  d\mu=\sum_{k=0}^{\infty}\sum
_{l=1}^{a_{k}}\int_{0}^{\infty}f_{k,l}\left(  r\right)  r^{-k}d\mu
_{k,l}\left(  r\right)  \label{eqneuTT}%
\end{equation}
holds for any continuous, polynomially bounded function $f:\mathbb{R}%
^{n}\rightarrow\mathbb{C}$; here $f_{k,l}\left(  r\right)  $ are the
Laplace-Fourier coefficients. Equation (\ref{eqneuTT}) will be the key for
defining \emph{polyharmonic Gau\ss--Jacobi cubatures} as we shall show below.

Another strong supporting evidence for the nice properties of the notion of
pseudo--positivity is the satisfactory solution of the question of determinacy
in Section \ref{Sdeterminacy}: we show that the representing measure $\mu$ of
a pseudo-positive definite functional $T:\mathbb{C}\left[  x_{1}%
,x_{2},...,x_{d}\right]  \rightarrow\mathbb{C}$ is unique in the class of all
pseudo-positive signed measures whenever each component functional $T_{k,l}$
defined in (\ref{defpseudo}) has a unique representing measure on $\left[
0,\infty\right)  $ in the sense of Stieltjes (for the precise definition see
Section \ref{Sdeterminacy}). And vice versa, if a pseudo--positive functional
$T$ is determinate in the class of all pseudo-positive signed measures and the
summability condition (\ref{maincond}) is satisfied, then each functional
$T_{k,l}$ is determinate in the sense of Stieltjes. The proof is essentially
based on the properties of the Nevanlinna extremal measures.

Let us illustrate the notion of pseudo-positivity in the case where the signed
measure $\mu$ has a continuous density $w\left(  x\right)  $ with respect to
the Lebesgue measure $dx$. We put $d\mu\left(  x\right)  =w\left(  x\right)
dx$ in (\ref{eqpoll}) and take into account that $Y_{k,l}\left(  x\right)
=\left|  x\right|  ^{k}Y_{k,l}\left(  \theta\right)  $, and $dx=r^{d-1}d\theta
dr,$ (for detailed computations see Proposition \ref{ThmLaplace}). We see that
the component functionals $T_{k,l}$ defined in (\ref{defpseudo}) are now given
by
\begin{equation}
T_{k,l}\left(  p\right)  =\int_{0}^{\infty}p\left(  r^{2}\right)
r^{k+d-1}w_{k,l}\left(  r\right)  dr, \label{eqneuneu2}%
\end{equation}
where $w_{k,l}\left(  r\right)  $ are the Laplace-Fourier coefficients of the
function $w$ as defined in (\ref{eqfourier}). From (\ref{eqneuneu2}) it is
obvious that the non--negativity of $w_{k,l}$ implies that the measure $\mu$
is pseudo--positive and the corresponding functional $T$ defined by
(\ref{eqpoll}) is pseudo--positive definite. We regard now $d\mu_{k,l}\left(
r\right)  =r^{k+d-1}w_{k,l}\left(  r\right)  dr$ as a univariate non-negative
measure which represents the functional $T_{k,l}.$

Now we are moving to our main theme, the construction of Gau\ss--Jacobi type
cubatures for pseudo-positive measures. Let us first recall some terminology:
By a \emph{cubature formula} one usually means a linear functional of the
form
\begin{equation}
C\left(  f\right)  :=\alpha_{1}f\left(  x_{1}\right)  +....+\alpha_{s}f\left(
x_{s}\right)  =\int\limits_{\mathbb{R}}(\sum_{j=1}^{s}\alpha_{j}\delta\left(
x-x_{j}\right)  )f\left(  x\right)  dx \label{eqcub}%
\end{equation}
defined on the set $C\left(  \mathbb{R}^{n}\right)  $, the set of all
continuous complex-valued functions on $\mathbb{R}^{n};$ here $\delta$ is the
Dirac delta function. The points $x_{1},...,x_{s}$ are called \emph{nodes} and
the coefficients $\alpha_{1},...,\alpha_{s}\in\mathbb{R}$ \emph{weights}. A
cubature formula $C\left(  \cdot\right)  $ is \emph{exact on a subspace }$U$
of $C\left(  \mathbb{R}^{n}\right)  $ with respect to a measure $\nu$ if
\begin{equation}
C\left(  f\right)  =\int f\left(  x\right)  d\nu\label{GJequality}%
\end{equation}
holds for all $f\in U.$ If $U_{s}$ is the set of all polynomials of degree
$\leq s,$ and the cubature is exact on $U_{s}$ but not on $U_{s+1},$ we say
that $C$ has \emph{order} $s.$ It is common to call a cubature in the case
$d=1$ a \emph{quadrature}.

In our construction we will use the Gau\ss-Jacobi quadrature, so let us recall
its definition: Let $\nu$ be a non--negative measure on the interval $\left[
0,R\right]  $ and $s\geq1$ be an integer. If the cardinality of the support of
$\nu$ is $>s$ then there exist $s$ different points $t_{j}$ in the interval
$\left(  0,R\right)  $ and $s$ positive weights $\alpha_{j}$ which define the
classical \emph{Gau\ss--Jacobi measure }
\begin{equation}
d\nu^{\left(  s\right)  }\left(  t\right)  =\sum_{j=1}^{s}\alpha_{j}%
\delta\left(  t-t_{j}\right)  dt \label{GJquadrature}%
\end{equation}
and the corresponding \emph{Gau\ss--Jacobi quadrature }$C\left(  f\right)
=\int_{0}^{R}f\left(  t\right)  d\nu^{\left(  s\right)  }\left(  t\right)
=\sum_{j=1}^{s}\alpha_{j}f\left(  t_{j}\right)  $ satisfies (\ref{GJequality})
for all polynomials of degree $\leq2s-1.$ For formal reasons we put
$\nu^{\left(  s\right)  }\equiv\nu$ if the cardinality of $\nu$ is $\leq
s$.\footnote{The points $t_{j}$ are the zeros of the polynomial $Q^{s}\left(
t\right)  $ which is the $s-$th orthogonal with respect to the measure $\nu$
on $\left[  0,R\right]  .$ It is important that $t_{j}\in\left(  0,R\right)  $
if the support of $\nu$ has cardinality $>s,$ cf. Chapter $1$, Theorem $5.2$
in \cite{chiharabook}, and Theorem $5.1$ in Chapter $3.5$ in \cite{KrNu77}.}

It is not our intention to survey the numerous approaches to cubature
formulas; one may consult the references in \cite{Myso81}, \cite{DuXu01},
\cite{Stroud}, \cite{DavisRabinowitz} and in particular the monograph of S. L.
Sobolev \cite{Sob} where the minimization of the error functional of the
formula in (\ref{GJequality}) is the main objective.

Our approach, which is rather different from the usual cubature formulas, is
based on the identity (\ref{eqneuTT}). In order to be precise, let us begin
with the assumptions: let $\mu$ be a pseudo--positive measure and define the
non-negative \emph{component measures} $\mu_{k,l}$ by the identity
\begin{equation}
\int_{0}^{\infty}h\left(  \left|  x\right|  \right)  d\mu_{k,l}=\int
_{\mathbb{R}^{n}}h\left(  \left|  x\right|  \right)  Y_{k,l}\left(  x\right)
d\mu\label{componentmeasure}%
\end{equation}
valid for all continuous functions $h:\left[  0,\infty\right)  \rightarrow
\mathbb{C}$ with compact support. Let $\psi:\left[  0,\infty\right)
\rightarrow\left[  0,\infty\right)  $ be the transformation $\psi\left(
t\right)  =t^{2}$ and let $\mu_{k,l}^{\psi}$ be the image measure of
$\mu_{k,l}$ under $\psi$ (see (\ref{imagemeasure1}) below). Then
(\ref{eqneuTT}) becomes
\begin{equation}
\int_{\mathbb{R}^{n}}f\left(  x\right)  d\mu=\sum_{k=0}^{\infty}\sum
_{l=1}^{a_{k}}\int_{0}^{\infty}f_{k,l}\left(  \sqrt{t}\right)  t^{-\frac{1}%
{2}k}d\mu_{k,l}^{\psi}\left(  r\right)  . \label{eqneuTT+}%
\end{equation}
The \emph{main idea} is simple and consists in replacing in formula
(\ref{eqneuTT+}) the non-negative univariate measures $\mu_{k,l}^{\psi}$ by
their univariate Gau\ss-Jacobi quadratures\footnote{The upper index $s$ will
indicate the cardinality of the support.} $\nu_{k,l}^{\left(  s\right)  }$ of
order $2s-1.$ Let $\psi^{-1}$ be the inverse map of $\psi$ and put
$\sigma_{k,l}^{\left(  s\right)  }=\left(  \nu_{k,l}^{\left(  s\right)
}\right)  ^{\psi^{-1}}.$ Then we obtain a pseudo-positive definite functional
$T^{\left(  s\right)  }$ by setting
\begin{align}
T^{\left(  s\right)  }\left(  f\right)   &  :=\sum_{k=0}^{\infty}\sum
_{l=1}^{a_{k}}\int_{0}^{\infty}f_{k,l}\left(  r\right)  r^{-k}d\sigma
_{k,l}^{\left(  s\right)  }\left(  r\right) \label{eqTsss}\\
&  =\sum_{k=0}^{\infty}\sum_{l=1}^{a_{k}}\int_{0}^{\infty}f_{k,l}\left(
\sqrt{t}\right)  t^{-\frac{1}{2}k}d\nu_{k,l}^{\left(  s\right)  }\left(
t\right)  .\nonumber
\end{align}
As we have made it clear above, for $f\in\mathbb{C}\left[  x_{1}%
,x_{2},...,x_{d}\right]  $ we have $f_{k,l}\left(  r\right)  r^{-k}%
=p_{k,l}\left(  r^{2}\right)  ,$ where $p_{k,l}$ are polynomials. Hence, the
functional $T^{\left(  s\right)  }$ is well-defined on $\mathbb{C}\left[
x_{1},x_{2},...,x_{d}\right]  ,$ since then the series is finite. In fact, the
most important thing is to find a condition on the measure $\mu$ which
provides convergence of the series in (\ref{eqTsss}) for the class of
continuous, polynomially bounded functions $f.$

Since this is a very central result of our paper, let us give the \emph{main
argument }for proving the convergence of the series in (\ref{eqTsss}) in the
important case when all measures $\mu_{k,l}$ have their supports in the
compact interval $\left[  0,R\right]  .$ For the Laplace-Fourier coefficient,
defined in (\ref{eqfourier}), we have the simple estimate
\[
\left|  f_{k,l}\left(  r\right)  \right|  \leq C\max_{\left|  x\right|  \leq
R}\left|  f\left(  x\right)  \right|  \text{ for }0\leq r\leq R,
\]
based on the Cauchy inequality and the orthonormality of $\left\{
Y_{k,l}\left(  \theta\right)  \right\}  .$ Hence,
\[
\left|  \int_{0}^{\infty}f_{k,l}\left(  r\right)  r^{-k}d\sigma_{k,l}^{\left(
s\right)  }\left(  r\right)  \right|  \leq C\max_{\left|  x\right|  \leq
R}\left|  f\left(  x\right)  \right|  \int_{0}^{\infty}r^{-k}d\sigma
_{k,l}^{\left(  s\right)  }\left(  r\right)
\]
and
\begin{equation}
\left|  T^{\left(  s\right)  }\left(  f\right)  \right|  \leq C\max_{\left|
x\right|  \leq R}\left|  f\left(  x\right)  \right|  \sum_{k=0}^{\infty}%
\sum_{l=1}^{a_{k}}\int_{0}^{\infty}r^{-k}d\sigma_{k,l}^{\left(  s\right)
}\left(  r\right)  . \label{Rieszinequality}%
\end{equation}
Now here is the crux of the whole matter: for the convergence in
(\ref{eqTsss}) it would suffice to prove the inequality
\begin{equation}
\int_{0}^{\infty}r^{-k}d\sigma_{k,l}^{\left(  s\right)  }\left(  r\right)
\leq\int_{0}^{\infty}r^{-k}d\mu_{k,l}\left(  r\right)  .
\label{Chebinequality}%
\end{equation}
The famous Chebyshev extremal property\footnote{This has been proved by A.
Markov \cite{markovActa} and T. Stieltjes, cf. \cite[Chapter 4]{KrNu77} and
\cite[Chapter 3]{karlinstudden}.} of the Gau\ss--Jacobi quadrature provides us
with a proof of (\ref{Chebinequality}). So we see that the convergence of the
series in (\ref{eqTsss}) is a consequence of the summability condition
(\ref{maincond}) with $N=0.$ Further note that (\ref{Rieszinequality}) shows
that $T^{\left(  s\right)  }$ is a continuous functional: by the Riesz
representation theorem we infer the existence of a signed measure
$\sigma^{\left(  s\right)  }$ with support in the closed ball $B_{R}:=\left\{
x\in\mathbb{R}^{n}:\left|  x\right|  \leq R\right\}  $ such that
\[
T^{\left(  s\right)  }\left(  f\right)  =\int_{B_{R}}f\left(  x\right)
d\sigma^{\left(  s\right)  }\left(  x\right)  \
\]
for all continuous functions $f:B_{R}\rightarrow\mathbb{C}.$ Moreover, the
component measures of the pseudo--positive measure $\sigma^{\left(  s\right)
} $ are exactly the univariate measures $\sigma_{k,l}^{\left(  s\right)  }.$
The precise result is contained in Theorem \ref{ThmMain} in Section
\ref{Scubature}. In the case when not all measures $\mu_{k,l}$ have their
supports in a compact interval $\left[  0,R\right]  $ the argumentation has to
be modified, and one needs (\ref{maincond}) for all $N\geq0.$ The details are
provided in Section \ref{Smomentproblem} and Section \ref{Scubature}.

The exactness of the Gau\ss-Jacobi quadratures $\nu_{k,l}^{\left(  s\right)
}$ for polynomials of degree $\leq2s-1$ implies that $T^{\left(  s\right)  }$
and $\mu$ coincide on the set of all polynomials $P$ such that $\Delta
^{2s}P=0.$ This is due to the fact that in the Laplace--Fourier expansion
(\ref{frtheta}) the coefficients are given by $f_{k,l}\left(  r\right)
=r^{k}p_{k,l}\left(  r^{2}\right)  $ where $p_{k,l}$ are polynomials of degree
$2s-1.$ For that reason we call the measure $\sigma^{\left(  s\right)  }$ the
\emph{polyharmonic Gau\ss--Jacobi measure }or the \emph{polyharmonic
Gau\ss--Jacobi cubature of order }$s.$

In the following we want to discuss the properties of the \emph{polyharmonic
Gau\ss--Jacobi cubature} and it is natural to compare them with those of the
univariate Gau\ss--Jacobi quadrature. Among the various existing quadratures
(e.g. Newton-Cotes quadratures), the Gau\ss--Jacobi quadrature has the eminent
property that the weights are positive. This in turn is the key to prove the
convergence of the quadrature (see e.g. the discussion in \cite[p. 353]{Davis}
based on the theorems of P\'{o}lya and Steklov).

\begin{property}
\label{GJ2} (Stieltjes) For every continuous function $f$ the Gau\ss--Jacobi
quadrature $\int_{a}^{b}fd\nu^{\left(  s\right)  }$ converges to $\int_{a}%
^{b}fd\nu$, when $s$ tends to infinity.
\end{property}

A second important property of the Gau\ss--Jacobi quadrature is the error
estimate due to A. Markov (see \cite[p. 344]{Davis})

\begin{property}
\label{GJ3}(Markov) Let $\nu$ be a non-negative measure on $\left[
a,b\right]  $ whose support has cardinality $>s$. Then for any $2s$-times
continuously differentiable function $f:\left[  a,b\right]  \rightarrow
\mathbb{R}$ there exists $\xi\in\left(  a,b\right)  $ such that
\begin{equation}
\int_{a}^{b}f\left(  t\right)  d\nu\left(  t\right)  -\int_{a}^{b}f\left(
t\right)  d\nu^{\left(  s\right)  }\left(  t\right)  =\frac{1}{\left(
2s\right)  !}f^{\left(  2s\right)  }\left(  \xi\right)  \int_{a}^{b}\left|
Q_{s}\left(  t\right)  \right|  ^{2} \label{mark}%
\end{equation}
where $Q_{s}$ is the $s$-th orthogonal polynomial with respect to $\nu,$ with
leading coefficient $1$.
\end{property}

It is an amazing and non-trivial fact that properties \ref{GJ2}) and
\ref{GJ3}) have analogs for the polyharmonic Gau\ss--Jacobi cubature although
the approximation measures $\sigma^{\left(  s\right)  }$ are in general signed
measures. In Theorem \ref{TStieltjes} we show that $C_{s}\left(  f\right)  $
converges to $\int fd\mu$ for every continuous function $f:\mathbb{R}%
^{n}\rightarrow\mathbb{C}$. This property implies the numerical stability of
our cubature formula. In Section \ref{Serror} we prove an estimate for the
difference
\[
\mu\left(  f\right)  -T^{\left(  s\right)  }\left(  f\right)
\]
for functions $f\in C^{2s}\left(  \mathbb{R}^{d}\right)  $ by their
derivatives in the ball $B_{R}$ based on Markov's error estimate.

Let us outline the structure of the paper: In Section \ref{Smomentproblem} we
introduce the notions of pseudo--positive definite functional and
pseudo--positive measure, and we prove basic results about them. In Section
\ref{Sdeterminacy} we consider the determinacy question. In Section
\ref{Scubature} the polyharmonic Gau\ss--Jacobi cubature formula is presented
in detail. Section \ref{Serror} is devoted to a multivariate generalization of
the Markov's error estimate for the polyharmonic Gau\ss-Jacobi cubature.

Section \ref{Sannulus} and \ref{Sstrip} are devoted to definition of
polyharmonic Gau\ss--Jacobi cubatures in other domains with symmetries as the
annulus and the cylinder (periodic strip). In Section \ref{Sannulus} we
construct a Gau\ss-Jacobi cubature for pseudo-positive measures with support
in a closed annulus $A_{\rho,R}$ which is exact on the space of all functions
continuous on the closed annulus $A_{\rho,R}$ and polyharmonic of order $2s$
in the interior. While the case of the annulus is somewhat similar to that of
the ball, we have to introduce a new notion of pseudo-positivity in the case
of the cylinder (periodic strip) in Section \ref{Sstrip} in order to obtain
cubatures which preserve polyharmonic functions of order $2s.$ Moreover it is
not possible to use in the proof the usual univariate Gau\ss--Jacobi
quadratures; instead we need the existence of quadratures of Gau\ss
--Jacobi-type for Chebyshev systems (Theorem \ref{TMarkov}). The analog to the
crucial inequality (\ref{Chebinequality}) follows from the Markov--Krein
theory of extremal problems for the moment problem for Chebyshev
systems.\footnote{We use the name ''Markov--Krein theory'' following
\cite[Chapter 3]{karlinstudden}, while in \cite{KrNu77} this is called
''Chebyshev--Markov problem''.}

In Section \ref{Sexamples} we give explicit examples illustrating our results
and provide miscellaneous properties of pseudo--positive measures. In the last
Section \ref{Sfinal} we discuss shortly aspects of numerical implementation
and some background information about the polyharmonicity concept.

Finally, let us introduce some notations: the space of all continuous
complex-valued functions on a topological space $X$ is denoted by $C\left(
X\right)  .$ By $C_{c}\left(  X\right)  $ we denote the set of all $f\in
C\left(  X\right)  $ having compact support. Further $C_{pol}\left(
\mathbb{R}^{d}\right)  $ is the space of all polynomially bounded, continuous
functions, so for each $f\in C_{pol}\left(  \mathbb{R}^{d}\right)  $ there
exists $N\in\mathbb{N}_{0},$ such that $\left|  f\left(  x\right)  \right|
\leq C_{N}\left(  1+\left|  x\right|  \right)  ^{N}$ for some constant $C_{N}$
(depending on $f$ ) for all $x\in\mathbb{R}^{d}.$ Further, we define an useful
space of test functions
\begin{equation}
C^{\times}\left(  \mathbb{R}^{d}\right)  :=\{\sum_{k=0}^{N}\sum_{l=1}^{a_{k}%
}f_{k,l}\left(  \left|  x\right|  \right)  Y_{k,l}\left(  x\right)
:N\in\mathbb{N}_{0}\text{ and }f_{k,l}\in C\left[  0,\infty\right)  \}.
\label{Ccross}%
\end{equation}
which can be rephrased as the set of all continuous functions with a finite
Laplace-Fourier series. Moreover we set
\begin{equation}
C_{c}^{\times}\left(  \mathbb{R}^{d}\right)  :=C^{\times}\left(
\mathbb{R}^{d}\right)  \cap C_{c}\left(  \mathbb{R}^{d}\right)  .
\label{Ccrosscomp}%
\end{equation}
We need some terminology from measure theory: a signed measure on
$\mathbb{R}^{d}$ is a set function on the Borel $\sigma$-algebra on
$\mathbb{R}^{d}$ which takes real values and is $\sigma$-additive. For the
standard terminology, as Radon measure, Borel $\sigma$-algebra, etc., we refer
to \cite{BCR84}. By the \emph{Jordan decomposition } \cite[p. 125]{Cohn}, a
signed measure $\mu$ is the difference of two non-negative finite measures,
say $\mu=\mu^{+}-\mu^{-}$ with the property that there exist a Borel set $A$
such that $\mu^{+}\left(  A\right)  =0$ and $\mu^{-}\left(  \mathbb{R}%
^{n}\setminus A\right)  =0.$ The \emph{variation} of $\mu$ is defined as
$\left|  \mu\right|  :=\mu^{+}+\mu^{-}.$ The signed measure $\mu$ is called
\emph{moment measure} if all polynomials are integrable with respect to
$\mu^{+}$ and $\mu^{-},$ which is equivalent to integrability with respect to
the total variation. The \emph{support of a non-negative measure} $\mu$ on
$\mathbb{R}^{d}$ is defined as the complement of the largest open set $U$ such
that $\mu\left(  U\right)  =0.$ In particular, the \emph{support of the zero
measure} is the \emph{empty set}. The \emph{support of a signed measure}
$\sigma$ is defined as the support of the total variation $\left|
\sigma\right|  =\sigma_{+}+\sigma_{-}$ (see \cite[p. 226]{Cohn}). Recall that
in general, the supports of $\sigma_{+}$ and $\sigma_{-}$ are not disjoint
(cf. exercise 2 in \cite[p. 231]{Cohn}). For a surjective measurable mapping
$\varphi:X\rightarrow Y$ and a measure $\nu$ on $X$ the \emph{image measure}
$\nu^{\varphi}$ on $Y$ is defined by
\begin{equation}
\nu^{\varphi}\left(  B\right)  :=\nu\left(  \varphi^{-1}B\right)
\label{imagemeasure1}%
\end{equation}
for all Borel subsets $B$ of $Y.$ The equality $\int_{X}g\left(
\varphi\left(  x\right)  \right)  d\nu\left(  x\right)  =\int_{Y}g\left(
y\right)  d\nu^{\varphi}\left(  y\right)  $ holds for all integrable functions
$g$. We use the notation $\omega_{d-1}$ for the surface area of the unit
sphere, so
\begin{equation}
\omega_{d-1}:=\int_{\mathbb{S}^{d-1}}1d\theta. \label{omegad-1}%
\end{equation}
By $B_{R}$ we denote the closed ball $\left\{  x\in\mathbb{R}^{n}:\left|
x\right|  \leq R\right\}  .$ For $0\leq\rho<R\leq\infty$ we define the closed
annulus by
\begin{equation}
A_{\rho,R}:=\left\{  x\in\mathbb{R}^{d}:\rho\leq\left|  x\right|  \leq
R\right\}  . \label{annulus}%
\end{equation}
For a subset $B$ of $\mathbb{R}^{n}$ the interior is denoted by $B^{\circ}.$
Further, we denote the closed and the open interval respectively by $\left[
a,b\right]  =\left\{  x\in\mathbb{R}:a\leq x\leq b\right\}  $ and $\left(
a,b\right)  =\left\{  x\in\mathbb{R}:a<x<b\right\}  $.

\section{The moment problem for pseudo-positive definite
sequences\label{Smomentproblem}}

Let $T:\mathbb{C}\left[  x_{1},x_{2},...,x_{d}\right]  \rightarrow\mathbb{C}$
be a linear functional. For any solid harmonic polynomial $Y_{k,l}\left(
x\right)  $ we define the \emph{component functional} $T_{k,l}$ by
\begin{equation}
T_{k,l}\left(  p\right)  :=T_{c}\left(  p(\left|  x\right|  ^{2}%
)Y_{k,l}\left(  x\right)  \right)  \qquad\text{ for every }p\in\mathbb{C}%
\left[  x_{1}\right]  . \label{defpseudo2}%
\end{equation}
Let us give the precise definition of \emph{pseudo-positive definiteness,
}which we already mentioned in the introduction:

\begin{definition}
\label{DTpositivedef} A sequence $c=\left\{  c_{a}\right\}  _{\alpha
\in\mathbb{N}_{0}^{d}}$ , or the associated functional $T_{c},$ is
\emph{pseudo-positive definite} if for every $k\in\mathbb{N}_{0}$ and
$l=1,...,a_{k}$ the sequences $\{c_{j}^{\left(  k,l\right)  }\}_{j\in
\mathbb{N}_{0}}$ and $\{c_{j+1}^{\left(  k,l\right)  }\}_{j\in\mathbb{N}_{0}}$
defined in (\ref{defpseudoc}) are \emph{positive definite}. Clearly this is
the same to say that $T_{k,l}\left(  p^{\ast}\left(  t\right)  p\left(
t\right)  \right)  \geq0$ and $T_{k,l}\left(  t\cdot p^{\ast}\left(  t\right)
p\left(  t\right)  \right)  \geq0$ for every $p\left(  t\right)  \in
\mathbb{C}\left[  x_{1}\right]  $.
\end{definition}

First we recall the following result which may be found e.g. in
\cite{Baouendi} or \cite{Sob}.

\begin{proposition}
The Laplace-Fourier coefficient $f_{k,l}$ of a polynomial $f$ given by
(\ref{eqfourier}) is of the form $f_{k,l}\left(  r\right)  =r^{k}%
p_{k,l}\left(  r^{2}\right)  $ where $p_{k,l}$ is a univariate polynomial.
Hence, the Laplace-Fourier series (\ref{frtheta}) is equal to
\begin{equation}
f\left(  x\right)  =\sum_{k=0}^{\deg f}\sum_{l=1}^{a_{k}}p_{k,l}(\left|
x\right|  ^{2})Y_{k,l}\left(  x\right)  . \label{gauss}%
\end{equation}
\end{proposition}

Equality (\ref{gauss}) is a reformulation of the Gau\ss\ decomposition of a
polynomial which we have provided in (\ref{Gaussdecomposition}).

The next two Propositions characterize pseudo-positive definite sequences:

\begin{proposition}
\label{PropStieltjes} Let $c=\left\{  c_{\alpha}\right\}  _{\alpha
\in\mathbb{N}_{0}^{d}}$ be a pseudo-positive definite sequence and $T_{c}$ its
associated functional. Then for each $k\in\mathbb{N}_{0},$ $l=1,...,a_{k},$
there exist non-negative measures $\sigma_{k,l}$ with support in $\left[
0,\infty\right)  $ such that
\begin{equation}
T_{c}\left(  f\right)  =\sum_{k=0}^{\deg f}\sum_{l=1}^{a_{k}}\int_{0}^{\infty
}f_{k,l}\left(  r\right)  r^{-k}d\sigma_{k,l}\left(  r\right)
\label{eqTnicerep}%
\end{equation}
holds for all $f\in\mathbb{C}\left[  x_{1},x_{2},...,x_{d}\right]  $ where
$f_{k,l}\left(  r\right)  $, $k\in\mathbb{N}_{0},$ $l=1,...,a_{k},$ are the
Laplace-Fourier coefficients of $f.$
\end{proposition}

\begin{proof}
By the definition of pseudo-positive definiteness, $T_{k,l}\left(  p^{\ast
}\left(  t\right)  p\left(  t\right)  \right)  \geq0$ and $T_{k,l}\left(
t\cdot p^{\ast}\left(  t\right)  p\left(  t\right)  \right)  \geq0$ for each
univariate polynomial $p\left(  t\right)  $ where the component functional
$T_{k,l}$ is defined in (\ref{defpseudo2}). By the solution of the Stieltjes
moment problem there exists a non-negative measure $\mu_{k,l}$ with support in
$\left[  0,\infty\right)  $ representing the functional $T_{k,l},$ i.e.
satisfying
\begin{equation}
T_{k,l}\left(  p\right)  =\int_{0}^{\infty}p\left(  t\right)  d\mu
_{k,l}\left(  t\right)  \qquad\text{for every }p\in\mathbb{C}\left[  t\right]
. \label{eqT1}%
\end{equation}
Let now $\varphi:\left[  0,\infty\right)  \rightarrow\left[  0,\infty\right)
$ be defined by $\varphi\left(  t\right)  =\sqrt{t}.$ Then we put
$\sigma_{k,l}:=\mu_{k,l}^{\varphi}$ where $\mu_{k,l}^{\varphi}$ is the image
measure defined in (\ref{imagemeasure1}). We obtain
\begin{equation}
\int_{0}^{\infty}h\left(  t\right)  d\mu_{k,l}\left(  t\right)  =\int
_{0}^{\infty}h\left(  r^{2}\right)  d\mu_{k,l}^{\varphi}\left(  r\right)  .
\label{eqT2}%
\end{equation}
Now use (\ref{gauss}), the linearity of $T$ and the definition of $T_{k,l}$ in
(\ref{defpseudo2}), and the equations (\ref{eqT1}) and (\ref{eqT2}) to obtain
\[
T_{c}\left(  f\right)  =\sum_{k=0}^{\deg f}\sum_{l=1}^{a_{k}}T_{k,l}\left(
p_{k,l}\right)  =\sum_{k=0}^{\deg f}\sum_{l=1}^{a_{k}}\int_{0}^{\infty}%
p_{k,l}\left(  r^{2}\right)  d\mu_{k,l}^{\varphi}\left(  r\right)  .
\]
Since $p_{k,l}\left(  r^{2}\right)  =r^{-k}f_{k,l}\left(  r\right)  $ the
claim (\ref{eqTnicerep}) follows from the last equation, which ends the proof.
\end{proof}

The next result shows that the converse of Proposition \ref{PropStieltjes} is
also true; not less important, it is a natural way of defining pseudo-positive
definite sequences.

\begin{proposition}
\label{PropTTT}Let $\sigma_{k,l},$ $k\in\mathbb{N}_{0},$ $l=1,...,a_{k},$ be
non-negative moment measures with support in $\left[  0,\infty\right)  .$ Then
the functional $T:\mathbb{C}\left[  x_{1},x_{2},...,x_{d}\right]
\rightarrow\mathbb{C}$ defined by
\begin{equation}
T\left(  f\right)  :=\sum_{k=0}^{\deg f}\sum_{l=1}^{a_{k}}\int_{0}^{\infty
}f_{k,l}\left(  r\right)  r^{-k}d\sigma_{k,l} \label{defTTT}%
\end{equation}
is pseudo-positive definite, where $f_{k,l}\left(  r\right)  $, $k\in
\mathbb{N}_{0},$ $l=1,...,a_{k},$ are the Laplace-Fourier coefficients of $f.$
\end{proposition}

\begin{proof}
Let us compute $T_{k,l}\left(  p\right)  $ where $p$ is a univariate
polynomial: by definition, $T_{k,l}\left(  p\right)  =T\left(  p(\left|
x\right|  ^{2})Y_{k,l}\left(  x\right)  \right)  $. The Laplace-Fourier series
of the function $x\mapsto\left|  x\right|  ^{2j}p(\left|  x\right|
^{2})Y_{k,l}\left(  x\right)  $ is equal to $r^{2j}p\left(  r^{2}\right)
r^{k}Y_{k,l}\left(  \theta\right)  $, hence
\[
T_{k,l}\left(  t^{j}p\left(  t\right)  \right)  =T\left(  \left|  x\right|
^{2j}p(\left|  x\right|  ^{2})Y_{k,l}\left(  x\right)  \right)  =\int
_{0}^{\infty}r^{j}p\left(  r^{2}\right)  d\sigma_{k,l}%
\]
for every natural number $j.$ Taking $j=0$ and $j=1$ one concludes that
$T_{k,l}\left(  p^{\ast}\left(  t\right)  p\left(  t\right)  \right)  \geq0$
and $T_{k,l}\left(  tp^{\ast}\left(  t\right)  p\left(  t\right)  \right)
\geq0$ for all univariate polynomials $p$, hence $T$ is pseudo-positive definite.
\end{proof}

The \emph{pseudo-positive definiteness} is defined for a \emph{functional} on
$\mathbb{C}\left[  x_{1},x_{2},...,x_{d}\right]  .$ Now we introduce the
concept of \emph{pseudo-positivity} of a \emph{measure}:

\begin{definition}
\label{Dmupsp}A signed measure $\mu$ on $\mathbb{R}^{n}$ is called
\textbf{pseudo-positive} if
\begin{equation}
\int_{\mathbb{R}^{d}}h\left(  \left|  x\right|  \right)  Y_{k,l}\left(
x\right)  d\mu\left(  x\right)  \geq0 \label{defpspos2}%
\end{equation}
holds for every non-negative continuous function $h:\left[  0,\infty\right)
\rightarrow\left[  0,\infty\right)  $ with compact support.
\end{definition}

At first we need some basic properties of pseudo-positive measures.

\begin{proposition}
\label{pseudopos}Let $\mu$ be a pseudo-positive moment measure on
$\mathbb{R}^{d}.$ Then there exist unique moment measures $\mu_{k,l}$ defined
on $\left[  0,\infty\right)  ,$ which we call \emph{component measures}, such
that
\begin{equation}
\int_{0}^{\infty}h\left(  t\right)  d\mu_{k,l}\left(  t\right)  =\int
_{\mathbb{R}^{d}}h\left(  \left|  x\right|  \right)  Y_{k,l}\left(  x\right)
d\mu\label{eqlim}%
\end{equation}
holds for all $h\in C_{pol}\left[  0,\infty\right)  $. Further for each $f\in
C_{pol}^{\times}\left(  \mathbb{R}^{d}\right)  $
\[
\int_{\mathbb{R}^{d}}f\left(  x\right)  d\mu=\sum_{k=0}^{\infty}\sum
_{l=1}^{a_{k}}\int_{0}^{\infty}f_{k,l}\left(  r\right)  r^{-k}d\mu_{k,l}.
\]
\end{proposition}

\begin{proof}
By definition of pseudo-positivity, $M_{k,l}\left(  h\right)  :=\int
_{\mathbb{R}^{d}}h\left(  \left|  x\right|  \right)  Y_{k,l}\left(  x\right)
d\mu$ defines a positive functional on $C_{c}\left(  \left[  0,\infty\right)
\right)  .$ By the Riesz representation theorem there exists a unique
non-negative measure $\mu_{k,l}$ such that $M_{k,l}\left(  h\right)  =\int
_{0}^{\infty}h\left(  t\right)  d\mu_{k,l}$ for all $h\in C_{c}\left(  \left[
0,\infty\right)  \right)  .$ We want to show that (\ref{eqlim}) holds for all
$h\in C_{pol}\left[  0,\infty\right)  $. For this, let $u_{R}:\left[
0,\infty\right)  \rightarrow\left[  0,1\right]  $ be a \emph{cut--off
function, }so $u_{R}$ is continuous and decreasing such that
\begin{equation}
u_{R}\left(  r\right)  =1\text{ for all }0\leq r\leq R\text{ and }u_{R}\left(
r\right)  =0\text{ for all }r\geq R+1.\label{defuR}%
\end{equation}
Let $h\in C_{pol}\left[  0,\infty\right)  .$ Then $u_{R}h\in C_{c}\left(
\left[  0,\infty\right)  \right)  $ and
\begin{equation}
\int_{0}^{\infty}u_{R}\left(  t\right)  h\left(  t\right)  d\mu_{k,l}%
=\int_{\mathbb{R}^{d}}u_{R}\left(  \left|  x\right|  \right)  h\left(  \left|
x\right|  \right)  Y_{k,l}\left(  x\right)  d\mu.\label{eqlimitR}%
\end{equation}
Note that $\left|  u_{R}\left(  t\right)  h\left(  t\right)  \right|
\leq\left|  u_{R+1}\left(  t\right)  h\left(  t\right)  \right|  $ for all
$t\in\left[  0,\infty\right)  .$ Hence by the monotone convergence theorem
\begin{equation}
\int_{0}^{\infty}\left|  h\left(  t\right)  \right|  d\mu_{k,l}=\lim
_{R\rightarrow\infty}\int_{0}^{\infty}\left|  u_{R}\left(  t\right)  h\left(
t\right)  \right|  d\mu_{k,l}.\label{eqlimitR2}%
\end{equation}
On the other hand, it is obvious that
\begin{equation}
\left|  \int_{\mathbb{R}^{d}}u_{R}\left(  \left|  x\right|  \right)  \left|
h\left(  \left|  x\right|  \right)  \right|  Y_{k,l}\left(  x\right)
d\mu\right|  \leq\int_{\mathbb{R}^{d}}\left|  h\left(  \left|  x\right|
\right)  Y_{k,l}\left(  x\right)  \right|  d\left|  \mu\right|  .\label{eq36}%
\end{equation}
The last expression is finite since $\mu$ is a moment measure. From
(\ref{eqlimitR2}), (\ref{eqlimitR}) applied to $\left|  h\right|  $ and
(\ref{eq36}) it follows that $\left|  h\right|  $ is integrable for $\mu
_{k,l}.$ Using Lebesgue's convergence theorem for $\mu$ and (\ref{eqlimitR})
it is easy to that (\ref{eqlim}) holds. For the last statement note that each
$f\in C_{pol}^{\times}\left(  \mathbb{R}^{d}\right)  $ has a finite
Laplace-Fourier series, and it is easy to see that the Laplace-Fourier
coefficients $f_{k,l}$ are in $C_{pol}\left[  0,\infty\right)  $, see
(\ref{eqffpol}) below.
\end{proof}

The next theorem is the main result of this section and it provides a simple
sufficient condition for the pseudo-positive definite functional on
$\mathbb{C}\left[  x_{1},x_{2},...,x_{d}\right]  $ defined in (\ref{defTTT})
to possess a pseudo--positive representing measure. Let us note that not every
pseudo-positive definite functional has a pseudo-positive representing
measure, see Section \ref{Sexamples} for an example.

\begin{theorem}
\label{ThmRepG}Let $\sigma_{k,l},$ $k\in\mathbb{N}_{0},$ $l=1,...,a_{k},$ be
non-negative measures with support in $\left[  0,\infty\right)  $ such that
for any $N\in\mathbb{N}_{0}$
\begin{equation}
C_{N}:=\sum_{k=0}^{\infty}\sum_{l=1}^{a_{k}}\int_{0}^{\infty}r^{N}%
r^{-k}d\sigma_{k,l}<\infty\text{ .} \label{neu21}%
\end{equation}
Then for the functional $T:\mathbb{C}\left[  x_{1},x_{2},...,x_{d}\right]
\rightarrow\mathbb{C}$ defined by (\ref{defTTT}) there exists a
pseudo-positive, signed moment measure $\sigma$ such that
\[
T\left(  f\right)  =\int_{\mathbb{R}^{n}}fd\sigma\text{ for all }%
f\in\mathbb{C}\left[  x_{1},x_{2},...,x_{d}\right]  .
\]
\end{theorem}

\begin{remark}
\label{Rrep} 1. If the measures $\sigma_{k,l}$ have supports in the compact
interval $\left[  \rho,R\right]  $ for all $k\in\mathbb{N}_{0},$
$l=1,...,a_{k},$ then the measure $\sigma$ in Theorem \ref{ThmRepG} has
support in the annulus $\left\{  x\in\mathbb{R}^{d}:\rho\leq\left|  x\right|
\leq R\right\}  .$

2. In the case of $R<\infty$ , it obviously suffices to assume that
$C_{0}<\infty$ instead of $C_{N}<\infty$ for all $N\in\mathbb{N}_{0}.$

3. The proof of Theorem \ref{ThmRepG} shows that $\sigma_{k,l}$ is equal to
the component measure induced by $\sigma$ with respect to the solid harmonic
$Y_{k,l}\left(  x\right)  .$
\end{remark}

\begin{proof}
1. We show at first that $T$ can be extended to a linear functional
$\widetilde{T}$ defined on $C_{pol}\left(  \mathbb{R}^{d}\right)  $ by the
formula
\begin{equation}
\widetilde{T}\left(  f\right)  :=\sum_{k=0}^{\infty}\sum_{l=1}^{a_{k}}\int
_{0}^{\infty}f_{k,l}\left(  r\right)  r^{-k}d\sigma_{k,l} \label{eqtschlange}%
\end{equation}
for $f\in C_{pol}\left(  \mathbb{R}^{d}\right)  ,$ where $f_{k,l}\left(
r\right)  $ are the Laplace-Fourier coefficients of $f$. Indeed, since $f\in
C_{pol}\left(  \mathbb{R}^{d}\right)  $ is of polynomial growth there exists
$C>0$ and $N\in\mathbb{N}$ such that $\left|  f\left(  x\right)  \right|  \leq
C(1+\left|  x\right|  ^{N})$. If follows from (\ref{eqfourier}) that
\begin{equation}
\left|  f_{k,l}\left(  r\right)  \right|  \leq C\left(  1+r^{N}\right)
\sqrt{\omega_{d-1}}\sqrt{\int_{\mathbb{S}^{d-1}}\left|  Y_{k,l}\left(
\theta\right)  \right|  ^{2}d\theta}=C\left(  1+r^{N}\right)  \sqrt
{\omega_{d-1}}, \label{eqffpol}%
\end{equation}
where we used the Cauchy-Schwarz inequality and the fact that $Y_{k,l}$ is
orthonormal. Hence,
\[
\int_{0}^{\infty}\left|  f_{k,l}\left(  r\right)  \right|  r^{-k}d\sigma
_{k,l}\leq\sqrt{\omega_{d-1}}C\int_{0}^{\infty}\left(  1+r^{N}\right)
r^{-k}d\sigma_{k,l}.
\]
By assumption (\ref{neu21}) the latter integral exists, so $f_{k,l}\left(
r\right)  r^{-k}$ is integrable with respect to $\sigma_{k,l}.$ By summing
over all $k,l$ we obtain by (\ref{neu21}) that
\[
\sum_{k=0}^{\infty}\sum_{l=1}^{a_{k}}\left|  \int_{0}^{\infty}f_{k,l}\left(
r\right)  r^{-k}d\sigma_{k,l}\right|  <\infty,
\]
which implies the convergence of the series in (\ref{eqtschlange}). It follows
that $\widetilde{T}$ is well-defined.

2. Let $T_{0}$ be the restriction of the functional $\widetilde{T}$ to the
space $C_{c}\left(  \mathbb{R}^{d}\right)  $. We will show that $T_{0}$ is
continuous. Let $f\in C_{c}\left(  \mathbb{R}^{d}\right)  $ and suppose that
$f $ has support in the annulus $\left\{  x\in\mathbb{R}^{d}:\rho\leq\left|
x\right|  \leq R\right\}  $ (for the case $\rho=0$ this is a ball). Then by a
similar technique as above
\[
\left|  f_{k,l}\left(  r\right)  \right|  \leq\sqrt{\omega_{d-1}}\max
_{\rho\leq\left|  x\right|  \leq R}\left|  f\left(  x\right)  \right|  .
\]
Using (\ref{eqtschlange}) one arrives at
\begin{equation}
\left|  T_{0}\left(  f\right)  \right|  \leq\max_{\rho\leq\left|  x\right|
\leq R}\left|  f\left(  x\right)  \right|  \sqrt{\omega_{d-1}}\sum
_{k=0}^{\infty}\sum_{l=1}^{a_{k}}\int_{\rho}^{R}r^{-k}d\sigma_{k,l}.
\label{eqkey23}%
\end{equation}

3. First consider the case that all measures $\sigma_{k,l}$ have supports in
the interval $\left[  \rho,R\right]  $ with $R<\infty$ (cf. Remark
\ref{Rrep}). Then (\ref{eqkey23}) and the Riesz representation theorem for
compact spaces yield a representing measure with support in the annulus
$\left\{  x\in\mathbb{R}^{d}:\rho\leq\left|  x\right|  \leq R\right\}  .$ The
pseudo--positivity of $\mu$ will be proved in item 5.) below.

In the general case, we apply the Riesz representation theorem given in
\cite[p. 41, Theorem 2.5]{BCR84}: there exists a unique signed measure
$\sigma$ such that
\[
T_{0}\left(  g\right)  =\int_{\mathbb{R}^{d}}gd\sigma\qquad\text{ for all
}g\in C_{c}\left(  \mathbb{R}^{d}\right)  .
\]

4. Next we will show that the polynomials are integrable with respect to the
variation of the representation measure $\sigma.$ Let $\sigma=\sigma
_{+}-\sigma_{-}$ be the Jordan decomposition of $\sigma$. Following the
techniques of Theorem 2.4 and Theorem 2.5 in \cite[p. 42]{BCR84}, we have the
equality
\begin{equation}
\int_{\mathbb{R}^{d}}g\left(  x\right)  d\sigma_{+}=\sup\left\{  T_{0}\left(
h\right)  :h\in C_{c}\left(  \mathbb{R}^{d}\right)  \text{ with }0\leq h\leq
g\right\}  \label{defTnull}%
\end{equation}
which holds for any non-negative function $g\in C_{c}\left(  \mathbb{R}%
^{d}\right)  .$ Let $u_{R}$ be the cut-off function defined in (\ref{defuR}).
We want to estimate $\int_{\mathbb{R}^{d}}g\left(  x\right)  d\sigma_{+}$ for
the function $g:=\left|  x\right|  ^{N}u_{R}(\left|  x\right|  ^{2}).$ In view
of (\ref{defTnull}), let $h\in C_{c}\left(  \mathbb{R}^{d}\right)  $ with
$0\leq h\left(  x\right)  \leq\left|  x\right|  ^{N}u_{R}(\left|  x\right|
^{2})$ for all $x\in\mathbb{R}^{d}.$ Then for the Laplace-Fourier coefficient
$h_{k,l}$ of $h$ we have the estimate
\[
\left|  h_{k,l}\left(  r\right)  \right|  \leq\sqrt{\int_{\mathbb{S}^{d-1}%
}\left|  h\left(  r\theta\right)  \right|  ^{2}d\theta}\sqrt{\int
_{\mathbb{S}^{d-1}}\left|  Y_{k,l}\left(  \theta\right)  \right|  ^{2}d\theta
}\leq r^{N}u_{R}\left(  r^{2}\right)  \sqrt{\omega_{d-1}}.
\]
According to (\ref{eqtschlange})
\[
T_{0}\left(  h\right)  \leq\left|  T_{0}\left(  h\right)  \right|  \leq
\sqrt{\omega_{d-1}}\sum_{k=0}^{\infty}\sum_{l=1}^{a_{k}}\int_{0}^{\infty}%
r^{N}r^{-k}d\sigma_{k,l}=:D_{N}.
\]
From (\ref{defTnull}) it follows that $\int_{\mathbb{R}^{d}}\left|  x\right|
^{N}u_{R}(\left|  x\right|  ^{2})d\sigma_{+}\leq D_{N}$ for all $R>0$ (note
that $D_{N}$ does not depend on $R$ ). By the monotone convergence theorem
(note that $u_{R}\left(  x\right)  \leq u_{R+1}\left(  x\right)  $ for all
$x\in\mathbb{R}^{d})$ we obtain
\[
\int_{\mathbb{R}^{d}}\left|  x\right|  ^{N}d\sigma_{+}=\lim_{R\rightarrow
\infty}\int_{\mathbb{R}^{d}}\left|  x\right|  ^{N}u_{R}(\left|  x\right|
^{2})d\sigma_{+}\leq D_{N}.
\]
Similarly one shows that $\int_{\mathbb{R}^{d}}\left|  x\right|  ^{N}%
d\sigma_{-}<\infty$ by considering the functional $S=-T_{0}$. It follows that
all polynomials are integrable with respect to $\sigma_{+}$ and $\sigma_{-}$.
Using similar arguments it is not difficult to see that for all $g\in
C^{\times}\left(  \mathbb{R}^{d}\right)  \cap C_{pol}\left(  \mathbb{R}%
^{d}\right)  $
\begin{equation}
\int_{\mathbb{R}^{d}}g\left(  x\right)  d\sigma=\widetilde{T}\left(  g\right)
.\text{ } \label{eqidentmt}%
\end{equation}

5. It remains to prove that $\sigma$ is pseudo-positive as given by Definition
\ref{Dmupsp}. Let $h\in C_{c}\left(  \left[  0,\infty\right)  \right)  $ be a
non-negative function. The Laplace-Fourier coefficients $f_{k^{\prime
},l^{\prime}}$ of $f\left(  x\right)  :=h\left(  \left|  x\right|  \right)
Y_{k,l}\left(  x\right)  $ are given by $f_{k^{\prime}l^{\prime}}\left(
r\right)  =\delta_{kk^{\prime}}\delta_{ll^{\prime}}h\left(  r\right)  r^{k}$
and by (\ref{eqidentmt}) it follows that
\[
\int_{\mathbb{R}^{d}}h\left(  \left|  x\right|  \right)  Y_{k,l}\left(
x\right)  d\sigma=\widetilde{T}\left(  f\right)  =\int_{0}^{\infty}%
f_{k,l}\left(  r\right)  r^{-k}d\sigma_{k,l}=\int_{0}^{\infty}h\left(
r\right)  d\sigma_{k,l}.
\]
Since $\sigma_{k,l}$ are non-negative measures, the last term is non-negative.
According to definition (\ref{defpspos2}), $\sigma$ is pseudo-positive. The
proof is complete.
\end{proof}

The following is a solution to the modified moment problem as explained in the
introduction. It is an immediate consequence of Theorem \ref{ThmRepG}.

\begin{corollary}
\label{CorMoment}Let $T:\mathbb{C}\left[  x_{1},x_{2},...,x_{d}\right]
\rightarrow\mathbb{C}$ be a pseudo-positive definite functional. Let
$\sigma_{k,l},k\in\mathbb{N}_{0},$ $l=1,...,a_{k},$ be the non-negative
measures with supports in $\left[  0,\infty\right)  $ representing the
functional $T$ as obtained in Proposition \ref{PropStieltjes}. If for any
$N\in\mathbb{N}_{0}$
\begin{equation}
\sum_{k=0}^{\infty}\sum_{l=1}^{a_{k}}\int_{0}^{\infty}r^{N}r^{-k}d\sigma
_{k,l}<\infty, \label{eqneuCN}%
\end{equation}
then there exists a pseudo-positive, signed moment measure $\sigma$ such that
\[
T\left(  f\right)  =\int fd\sigma\qquad\text{ for all }f\in\mathbb{C}\left[
x_{1},x_{2},...,x_{d}\right]  .
\]
\end{corollary}

It would be interesting to see whether the summability condition
(\ref{eqneuCN}) may be weakened, cf. also the discussion at the end of Section
\ref{Sexamples}.

By the uniqueness of the representing measure in the Riesz representation
theorem for compact spaces we conclude from Theorem \ref{ThmRepG}:

\begin{corollary}
\label{ppp}Let $\mu$ be a signed measure with compact support. Then $\mu$ is
pseudo-positive if and only if $\mu$ is pseudo-positive definite as a
functional on $\mathbb{C}\left[  x_{1},x_{2},...,x_{d}\right]  .$
\end{corollary}

Let us remark that Corollary \ref{ppp} does not hold without the compactness
assumption which follows from well known arguments in the univariate case:
Indeed, let $\nu_{1}$ be a non-negative moment measure on $\left[
0,\infty\right)  $ which is not determined in the sense of Stieltjes; hence
there exists a non-negative moment measure $\nu_{2}$ on $\left[
0,\infty\right)  $ such that $\nu_{1}\left(  p\right)  =\nu_{2}\left(
p\right)  $ for all univariate polynomials. Since $\nu_{1}\neq\nu_{2}$ there
exists a continuous function $h:\left[  0,\infty\right)  \rightarrow\left[
0,\infty\right)  $ with compact support that $\nu_{1}\left(  h\right)  \neq
\nu_{2}\left(  h\right)  .$ Without loss of generality assume that
\begin{equation}
\int_{0}^{\infty}h\left(  r\right)  d\nu_{1}-\int_{0}^{\infty}h\left(
r\right)  d\nu_{2}<0. \label{eqverschieden}%
\end{equation}
For $i=1,2$ define $\mu_{i}=d\theta d\nu_{i},$ so for any $f\in C\left(
\mathbb{R}^{d}\right)  $ of polynomial growth
\[
\int fd\mu_{i}=\int_{0}^{\infty}\int_{\mathbb{S}^{d-1}}f\left(  r\theta
\right)  d\theta d\nu_{i}.
\]
For a polynomial $f$ let $f_{0}$ be the first Laplace--Fourier coefficient.
Then $\int fd\mu_{i}=\int_{0}^{\infty}f_{0}\left(  r\right)  d\nu_{i}$ for
$i=1,2.$ Since $\nu_{1}\left(  p\right)  =\nu_{2}\left(  p\right)  $ for all
univariate polynomials it follows that $\int fd\mu_{1}=\int fd\mu_{2}$ for all
polynomials. Then $\mu:=\mu_{1}-\mu_{2}$ is a signed measure which is
pseudo-positive definite since $\mu\left(  P\right)  =0$ for all polynomials
$P.$ It is not pseudo-positive since $\mu_{0}\left(  h\right)  =\int h\left(
\left|  x\right|  \right)  d\mu<0$ by (\ref{eqverschieden}).

\section{Determinacy for pseudo-positive definite
functionals\label{Sdeterminacy}}

Let $M^{\ast}\left(  \mathbb{R}^{d}\right)  $ be the set of all \emph{signed
moment measures}, and $M_{+}^{\ast}\left(  \mathbb{R}^{d}\right)  $ be the set
of \emph{non--negative moment measures} on $\mathbb{R}^{d}$. On $M^{\ast
}\left(  \mathbb{R}^{d}\right)  $ we define an equivalence relation: we say
that $\sigma\sim\mu$ \ for two elements $\sigma,\mu\in M^{\ast}\left(
\mathbb{R}^{d}\right)  $ if and only if $\int_{\mathbb{R}^{d}}fd\sigma
=\int_{\mathbb{R}^{d}}fd\mu$ for all $f\in\mathbb{C}\left[  x_{1}%
,x_{2},...,x_{d}\right]  .$

\begin{definition}
\label{Ddetermined}Let $\mu\in$ $M^{\ast}\left(  \mathbb{R}^{d}\right)  $ be a
pseudo-positive measure. We define
\[
V_{\mu}=\left\{  \sigma\in M^{\ast}\left(  \mathbb{R}^{d}\right)
:\sigma\text{ is pseudo-positive and }\sigma\sim\mu\right\}  .
\]
We say that the measure $\mu\in$ $M^{\ast}\left(  \mathbb{R}^{d}\right)  $ is
\emph{determined in the class of pseudo-positive measures} if $V_{\mu}$ has
only one element, i.e. is equal to $\left\{  \mu\right\}  .$
\end{definition}

Recall that a positive definite functional $\phi:\mathcal{P}_{1}%
\rightarrow\mathbb{R}$ is \emph{determined in the sense of Stieltjes} if the
set
\begin{equation}
W_{\phi}^{Sti}:=\left\{  \tau\in M_{+}^{\ast}\left(  \left[  0,\infty\right)
\right)  :\int_{0}^{\infty}r^{m}d\tau=\phi\left(  r^{m}\right)  \text{ for all
}m\in\mathbb{N}_{0}\right\}  \label{defES}%
\end{equation}
has exactly one element, cf. \cite[p. 210]{BeTh91}.

According to Proposition \ref{pseudopos}, we can associate to a
pseudo-positive measure $\mu$ the sequence of non-negative component measures
$\mu_{k,l},k\in\mathbb{N}_{0},l=1,..,a_{k}$ with support in $\left[
0,\infty\right)  .$ The measures $\mu_{k,l}$ contain all information about
$\mu.$ Indeed, we prove

\begin{proposition}
\label{soso}Let $\mu$ and $\sigma$ be pseudo-positive measures and let
$\mu_{k,l}$ and $\sigma_{k,l}$ be as in Proposition \ref{pseudopos}. If
$\mu_{k,l}=\sigma_{k,l}$ for all $k\in\mathbb{N}_{0},l=1,..,a_{k}$ then
$\mu=\sigma.$
\end{proposition}

\begin{proof}
Let $h\in C_{c}\left[  0,\infty\right)  .$ Then, using the assumption
$\mu_{k,l}=\sigma_{k,l},$ we obtain
\[
\int_{\mathbb{R}^{d}}h\left(  \left|  x\right|  \right)  Y_{k,l}\left(
x\right)  d\mu=\int_{0}^{\infty}h\left(  t\right)  d\mu_{k,l}=\int
_{\mathbb{R}^{d}}h\left(  \left|  x\right|  \right)  Y_{k,l}\left(  x\right)
d\sigma.
\]
Since each $f\in C_{c}^{\times}\left(  \mathbb{R}^{d}\right)  $ is a finite
linear combination of functions of the type $h\left(  \left|  x\right|
\right)  Y_{k,l}\left(  x\right)  $, we obtain that $\int_{\mathbb{R}^{d}%
}fd\mu=\int_{\mathbb{R}^{d}}fd\sigma$ for all $f\in C_{c}^{\times}\left(
\mathbb{R}^{d}\right)  .$ We apply Proposition \ref{Ldensity} to see that
$\mu$ is equal to $\sigma.$
\end{proof}

The following result is proved in \cite[Proposition 3.1]{BeTh91}:

\begin{proposition}
\label{Ldensity} Let $\mu$ and $\sigma$ be signed measures on $\mathbb{R}%
^{d}.$ If $\int_{\mathbb{R}^{d}}fd\mu=\int_{\mathbb{R}^{d}}fd\sigma$ for all
$f\in C_{c}^{\times}\left(  \mathbb{R}^{d}\right)  ,$ then $\mu$ is equal to
$\sigma.$
\end{proposition}

We can characterize $V_{\mu}$ in the case that only finitely many $\mu_{k,l}$
are nonzero.

\begin{theorem}
Let $\mu$ be a pseudo-positive measure on $\mathbb{R}^{n}$ such that
$\mu_{k,l}=0$ for all $k>k_{0},l=1,...,a_{k}.$ Then $V_{\mu}$ is affinely
isomorphic to the set
\begin{equation}
\oplus_{k=0}^{k_{0}}\oplus_{l=1}^{a_{k}}\{\rho_{k,l}\in W_{\mu_{k,l}^{\psi}%
}^{Sti}:\int_{0}^{\infty}t^{-\frac{1}{2}k}d\rho_{k,l}<\infty\}
\label{descript}%
\end{equation}
where the isomorphism is given by $\sigma\longmapsto\left(  \sigma_{k,l}%
^{\psi}\right)  _{k=1,..,k_{0},l=1,...,a_{k}}$ and the map $\psi:\left[
0,\infty\right)  \rightarrow\left[  0,\infty\right)  $ is defined by
$\psi\left(  t\right)  =t^{2},$ cf. (\ref{imagemeasure1}).
\end{theorem}

\begin{proof}
Let $\sigma$ be in $V_{\mu}.$ Let $\sigma_{k,l}$ and $\mu_{k,l}$ be the unique
moment measures obtained in Proposition \ref{pseudopos}. Then
\[
\int_{0}^{\infty}h\left(  t\right)  d\sigma_{k,l}^{\psi}=\int_{0}^{\infty
}h\left(  t^{2}\right)  d\sigma_{k,l}=\int_{\mathbb{R}^{n}}h(\left|  x\right|
^{2})Y_{k,l}\left(  x\right)  d\sigma\left(  x\right)
\]
for all $h\in C_{pol}\left[  0,\infty\right)  ,$ and an analog equation is
valid for $\mu_{k,l}$ and $\mu.$ Taking polynomials $h\left(  t\right)  $ we
see that $\sigma_{k,l}\in W_{\mu_{k,l}^{\psi}}^{Sti}$ using the assumption
that $\mu\sim\sigma.$ Using a simple approximation argument it is easy to see
from (\ref{eqlim}) that
\[
\int_{0}^{\infty}t^{-\frac{1}{2}k}d\sigma_{k,l}^{\psi}=\int_{\mathbb{R}^{n}%
}Y_{k,l}\left(  \frac{x}{\left|  x\right|  }\right)  d\sigma\left(  x\right)
.
\]
Since $x\longmapsto Y_{k,l}\left(  \frac{x}{\left|  x\right|  }\right)  $ is
bounded on $\mathbb{R}^{n}$, say by $M,$ we obtain the estimate
\[
\left|  \int_{0}^{\infty}t^{-\frac{1}{2}k}d\sigma_{k,l}^{\psi}\right|  \leq
M\int_{\mathbb{R}^{n}}1d\left|  \sigma\right|  <\infty.
\]
It follows that $\left(  \sigma_{k,l}^{\psi}\right)  _{k=1,..,k_{0}%
,l=1,...,a_{k}}$ is contained in the set on the right hand side in
(\ref{descript}).

Let now $\rho_{k,l}\in W_{\mu_{k,l}^{\psi}}^{Sti}$ be given such that
$\int_{0}^{\infty}t^{-\frac{1}{2}k}d\rho_{k,l}<\infty$ for $k=1,..,k_{0}%
,l=1,...,a_{k}.$ Define $\sigma_{k,l}=\rho_{k,l}^{\psi^{-1}}$ and
$\sigma_{k,l}=0$ for $k>k_{0}$. Then by Theorem \ref{ThmRepG} there exists a
measure $\tau\in V_{\mu}$ such that $\tau_{k,l}=\sigma_{k,l}.$ This shows the
surjectivity of the map. Let now $\sigma$ and $\tau$ are in $V_{\mu}$ with
$\sigma_{k,l}^{\psi}=\tau_{k,l}^{\psi}$ for $k=1,..,k_{0},l=1,...,a_{k}.$ The
property $\sigma\in V_{\mu}$ implies that $\sigma_{k,l}^{\psi}\in W_{\mu
_{k,l}^{\psi}}^{Sti}$ for all $k\in\mathbb{N}_{0},l=1,...,a_{k},$ hence
$\sigma_{k,l}^{\psi}=0$ for $k>k_{0},$ and similarly $\tau_{k,l}^{\psi}=0.$
Hence $\sigma_{k,l}=\tau_{k,l}$ for all $k\in\mathbb{N}_{0},l=1,...,a_{k},$
and this implies that $\sigma=\tau$ by Proposition \ref{soso}.
\end{proof}

The following is a sufficient condition for a functional $T$ to be determined
in the class of pseudo-positive measures.

\begin{theorem}
\label{ThmUniq}Let $T:\mathbb{C}\left[  x_{1},x_{2},...,x_{d}\right]
\rightarrow\mathbb{R}$ be a pseudo-positive definite functional. If the
functionals $T_{k,l}:\mathbb{C}\left[  x_{1}\right]  \rightarrow\mathbb{C}$
are determined in the sense of Stieltjes then there exists at most one
pseudo-positive, signed moment measure $\mu$ on $\mathbb{R}^{d}$ with
\begin{equation}
T\left(  f\right)  =\int_{\mathbb{R}^{d}}fd\mu\qquad\text{ for all }%
f\in\mathbb{C}\left[  x_{1},x_{2},...,x_{d}\right]  . \label{eqTident}%
\end{equation}
\end{theorem}

\begin{proof}
Let us suppose that $\mu$ and $\sigma$ are pseudo-positive, signed moment
measures on $\mathbb{R}^{d}$ representing $T.$ Taking $f=\left|  x\right|
^{2N}Y_{k,l}\left(  x\right)  $ we obtain from (\ref{eqTident}) that
\[
\int_{\mathbb{R}^{d}}\left|  x\right|  ^{2N}Y_{k,l}\left(  x\right)
d\mu=T_{k,l}\left(  t^{N}\right)  =\int_{\mathbb{R}^{d}}\left|  x\right|
^{2N}Y_{k,l}\left(  x\right)  d\sigma.
\]
for all $N\in\mathbb{N}_{0}.$ Let $\mu_{k,l}$ and $\sigma_{k,l}$ as in
Proposition \ref{pseudopos}, and consider $\psi:\left[  0,\infty\right)
\rightarrow\left[  0,\infty\right)  $ defined by $\psi\left(  t\right)
=t^{2}$. Then the image measures $\mu_{k,l}^{\psi}$ and $\sigma_{k,l}^{\psi}$
are non-negative measures with supports on $\left[  0,\infty\right)  $ such
that $\int_{0}^{\infty}t^{N}d\mu_{k,l}^{\psi}=T_{k,l}\left(  t^{N}\right)
=\int_{0}^{\infty}t^{N}d\sigma_{k,l}^{\psi}.$ Our assumption implies that
$\mu_{k,l}^{\psi}=\sigma_{k,l}^{\psi}$, so $\mu_{k,l}=\sigma_{k,l}$.
Proposition \ref{soso} implies that $\mu$ is equal to $\sigma$.
\end{proof}

In the following we want to prove the converse of the last theorem, which is
more subtle. We need now some special results about \emph{Nevanlinna extremal
measures}. Let us introduce the following notation: for a non-negative measure
$\phi\in M_{+}^{\ast}\left(  \mathbb{R}\right)  $ we put\footnote{Here in
order to avoid mixing of the notations, we retain the notation $\left[
\phi\right]  $ from the one--dimensional case in \cite{BeTh91}.}
\[
\left[  \phi\right]  :=\left\{  \sigma\in M_{+}^{\ast}\left(  \mathbb{R}%
\right)  :\text{ }\sigma\sim\phi\right\}  .
\]

\begin{proposition}
\label{PropSD} Let $\nu$ be a non-negative moment measure on $\mathbb{R}$ with
support in $\left[  0,\infty\right)  $ which is not determined in the sense of
Stieltjes, or applying the notation (\ref{defES}) $W_{\nu}^{Sti}\neq\left\{
\nu\right\}  .$ Then there exist uncountably many $\sigma\in W_{\nu}^{Sti}$
such that $\int_{0}^{\infty}u^{-k}d\sigma<\infty$ for all $k\in\mathbb{N}_{0}.$
\end{proposition}%

\proof
In the proof we will borrow some arguments about the Stieltjes problem as
given in \cite{Chih82} or \cite{Pede95}. As in the proof of Proposition 4.1 in
\cite{Pede95} let $\varphi:\left(  -\infty,\infty\right)  \rightarrow\left[
0,\infty\right)  $ be defined by $\varphi\left(  x\right)  =x^{2}.$ If
$\lambda$ is a measure on $\mathbb{R}$ define a measure $\lambda^{-}$ by
$\lambda^{-}\left(  A\right)  :=\lambda\left(  -A\right)  $ for each Borel set
$A$ where $-A:=\left\{  -x:x\in A\right\}  .$ The measure is symmetric if
$\lambda^{-}=\lambda.$ For each $\tau\in W_{\nu}^{Sti}$ define a measure
$\widetilde{\tau}:=\frac{1}{2}\left(  \tau^{\varphi}+\left(  \tau^{\varphi
}\right)  ^{-}\right)  $ which is clearly symmetric, in particular
$\widetilde{\nu}$ is symmetric. As pointed out in \cite{Pede95}, the map
$\widetilde{\cdot}:W_{\nu}^{Sti}\rightarrow\left[  \widetilde{\nu}\right]  \ $
is injective and the image is exactly the set of all symmetric measures in the
set $\left[  \widetilde{\nu}\right]  .$ The inverse map of $\widetilde{\cdot}$
defined on the image space is just the map $\sigma\rightarrow\sigma^{\varphi}.$

It follows that $\widetilde{v}$ is not determined, so we can make use of the
Nevanlinna theory for the indeterminate measure $\widetilde{\nu},$ see p. 54
in \cite{Akhi65}. We know by formula II.4.2 (9) and II.4.2 (10) in
\cite{Akhi65} that for every $t\in\mathbb{R}$ there exists a unique
Nevanlinna--extremal measure $\sigma_{t}$ such that
\[
\int_{-\infty}^{\infty}\frac{d\sigma_{t}\left(  u\right)  }{u-z}%
=-\frac{A\left(  z\right)  t-C\left(  z\right)  }{B\left(  z\right)
t-D\left(  z\right)  },
\]
where $A\left(  z\right)  ,B\left(  z\right)  ,C\left(  z\right)  ,D\left(
z\right)  $ are entire functions. Since the support of $\sigma_{t}$ is the
zero-set of the entire function $B\left(  z\right)  t-D\left(  z\right)  $ it
follows that the measure $\sigma_{t}$ has no mass in $0$ for $t\neq0,$ and now
it is clear that $\sigma_{t}([-\delta,\delta])=0$ for $t\neq0$ and suitable
$\delta>0$ (this fact is pointed out at least in the reference \cite[p.
210]{BeTh91}). It follows that
\begin{equation}
\int_{-\infty}^{\infty}\left|  u\right|  ^{-k}d\sigma_{t}<\infty
\label{eqbounded}%
\end{equation}
since the function $u\longmapsto\left|  u\right|  ^{-k}$ is bounded on
$\mathbb{R}\setminus\left[  -\delta,\delta\right]  $ for each $\delta>0$.
Using the fact that the functions $A\left(  z\right)  $ and $B\left(
z\right)  $ of the Nevanlinna matrix are odd, while the functions $B\left(
z\right)  $ and $C\left(  z\right)  $ are even, one derives that the measure
$\rho_{t}:=\frac{1}{2}\sigma_{t}+\frac{1}{2}\sigma_{-t}$ is symmetric. Further
from the equation $A\left(  z\right)  D\left(  z\right)  -B\left(  z\right)
C\left(  z\right)  =1$ it follows that $\rho_{t}\neq\rho_{s}$ for positive
numbers $t\neq s.$ By the above we know that $\rho_{t}^{\varphi}\neq\rho
_{s}^{\varphi}$. This finishes the proof.
\endproof

\begin{theorem}
Let $\mu$ be a pseudo-positive signed measure on $\mathbb{R}^{d}$ such that
the summability assumption (\ref{maincond}) holds. Then $V_{\mu}$ contains
exactly one element if and only if each $\mu_{k,l}^{\psi}$ is determined in
the sense of Stieltjes.
\end{theorem}

\begin{proof}
Let $\mu_{k,l}$ be the component measures as defined in Proposition
\ref{pseudopos}. Assume that $V_{\mu}=\left\{  \mu\right\}  $ but that some
$\tau:=\mu_{k_{0},l_{0}}^{\psi}$ is not determined in the sense of Stieltjes
where $\psi\left(  t\right)  =t^{2}$ for $t\in\left[  0,\infty\right)  .$ By
Proposition \ref{PropSD} there exists a measure $\sigma\in W_{\tau}^{Sti}$
such $\sigma\neq\tau$ and $\int_{0}^{\infty}r^{-k}d\sigma<\infty.$ By Theorem
\ref{ThmRepG} there exists a pseudo-positive moment measure $\widetilde{\mu}$
representing the functional
\[
\widetilde{T}\left(  f\right)  :=\sum_{k=0,k\neq k_{0}}^{\infty}%
\sum_{l=1,l\neq l_{0}}^{a_{k}}\int_{0}^{\infty}f_{k,l}\left(  r\right)
r^{-k}d\mu_{k,l}+\int_{0}^{\infty}f_{k_{0},l_{0}}\left(  r\right)
r^{-k}d\sigma^{\psi^{-1}}.
\]
Then $\widetilde{\mu}$ is different from $\mu$ since $\sigma^{\psi^{-1}}%
\neq\mu_{k_{0},l_{0}}$ and $\widetilde{\mu}\in V_{\mu}$ since $\sigma\in
W_{\tau}^{Sti}.$ This contradiction shows that $\mu_{k_{0},l_{0}}^{\psi}$ is
determined in the sense of Stieltjes. The sufficiency follows from Theorem
\ref{ThmUniq}. The proof is complete.
\end{proof}

\section{Polyharmonic Gau\ss--Jacobi cubatures \label{Scubature}}

In this section we will prove the main result of the paper, the existence of
the polyharmonic Gau\ss-Jacobi cubature of order $s.$ The proof is based on
application of the famous Chebyshev extremal property of the Gau\ss--Jacobi measure.

\begin{theorem}
\label{ThmMain}Let $0\leq\rho<R\leq\infty.$ Let $\mu$ be a pseudo-positive
signed measure with support in the closed annulus $A_{\rho,R}$ such that
\begin{equation}
\sum_{k=0}^{\infty}\sum_{l=1}^{a_{k}}\int_{\mathbb{R}^{d}}Y_{k,l}\left(
\frac{x}{\left|  x\right|  }\right)  d\mu<\infty. \label{condmu}%
\end{equation}
Then for each natural number $s$ there exists a unique pseudo-positive, signed
measure $\sigma^{\left(  s\right)  }$ with support in $A_{\rho,R}$ such that

(i) The support of each component measure $\sigma_{k,l}^{\left(  s\right)  }$
of $\sigma^{\left(  s\right)  }$ (defined by (\ref{eqlim}) ) has cardinality
$\leq s.$

(ii) $\int Pd\mu=\int Pd\sigma^{\left(  s\right)  }$ for all polynomials $P$
with $\Delta^{2s}P=0.$
\end{theorem}

\begin{proof}
By Proposition \ref{pseudopos} the following identity holds
\begin{equation}
\int_{\mathbb{R}^{d}}f\left(  x\right)  d\mu\left(  x\right)  =\sum
_{k=0}^{\infty}\sum_{l=1}^{a_{k}}\int_{0}^{\infty}f_{k,l}\left(  r\right)
r^{-k}d\mu_{k,l}\left(  r\right)  \label{defmmm}%
\end{equation}
for any $f\in C_{pol}^{\times}\left(  \mathbb{R}^{d}\right)  $ where
$\mu_{k,l}\left(  h\right)  =\int h\left(  \left|  x\right|  \right)
Y_{k,l}\left(  x\right)  d\mu\left(  x\right)  .$ It is clear that $\mu_{k,l}$
has support in the interval $\left[  \rho,R\right]  .$ If the cardinality of
the support of $\mu_{k,l}$ is $\leq s$ we define $\sigma_{k,l}^{\left(
s\right)  }:=\mu_{k,l}. $ If the cardinality is strictly larger than $s$ we
define $\sigma_{k,l}^{\left(  s\right)  }$ as the non-negative measure such
that
\begin{equation}
\int_{\rho}^{R}r^{2j}d\sigma_{k,l}^{\left(  s\right)  }\left(  r\right)
=\int_{\rho}^{R}r^{2j}d\mu_{k,l}\left(  r\right)  \label{eggauss}%
\end{equation}
for all $j=0,...,2s-1.$ The existence of $\sigma_{k,l}^{\left(  s\right)  }$
is proved as follows: Let $\psi:\left[  \rho,R\right]  \rightarrow\left[
\rho^{2},R^{2}\right]  $ be the map $\psi\left(  t\right)  =t^{2}$. Then the
image measure $\mu_{k,l}^{\psi}$ is a measure on $\left[  \rho^{2}%
,R^{2}\right]  $ and its support has clearly cardinality $>s.$ Let $\nu
_{k,l}^{\left(  s\right)  }$ be the Gau\ss--Jacobi quadrature of $\mu
_{k,l}^{\psi}.$ From the Gau\ss--Jacobi quadrature formula (\ref{GJquadrature}%
) (see also the footnote after it) follows that $\nu_{k,l}^{\left(  s\right)
}$ has support in the open interval $\left(  \rho^{2},R^{2}\right)  $ and
\[
\int_{\rho^{2}}^{R^{2}}t^{j}d\nu_{k,l}^{\left(  s\right)  }=\int_{\rho^{2}%
}^{R^{2}}t^{j}d\mu_{k,l}^{\psi}%
\]
for $j=0,...,2s-1$. Now it is easily seen that $\sigma_{k,l}^{\left(
s\right)  }:=\left(  \nu_{k,l}^{\left(  s\right)  }\right)  ^{\psi^{-1}}$
satisfies (\ref{eggauss}).

We now define a functional $T:\mathbb{C}\left[  x_{1},x_{2},...,x_{d}\right]
\rightarrow\mathbb{C}$ by putting
\begin{equation}
T^{\left(  s\right)  }\left(  f\right)  =\sum_{k=0}^{\deg f}\sum_{l=1}^{a_{k}%
}\int_{\rho}^{R}f_{k,l}\left(  r\right)  r^{-k}d\sigma_{k,l}^{\left(
s\right)  }. \label{defTss}%
\end{equation}

Let us show that
\begin{equation}
\int_{\mathbb{R}^{d}}P\left(  x\right)  d\mu\left(  x\right)  =T^{\left(
s\right)  }\left(  P\right)  \label{Tsp}%
\end{equation}
for all polynomials $P$ with $\Delta^{2s}P\left(  x\right)  =0.$ Indeed,
according to (\ref{gauss}) the Laplace-Fourier series of a polyharmonic
polynomial of order $2s$ can be written as
\begin{equation}
P\left(  x\right)  =\sum_{k=0}^{\deg f}\sum_{l=1}^{a_{k}}p_{k,l}\left(
r^{2}\right)  r^{k}Y_{k,l}\left(  \theta\right)  \label{eqP45}%
\end{equation}
and the univariate polynomials $p_{k,l}\left(  t\right)  $ have degree
$\leq2s-1$ (see e.g. \cite[Theorem 10.42, Remark 10.43]{Koun00}). Combining
(\ref{eqP45}) with (\ref{defTss}), (\ref{eggauss}) and (\ref{defmmm}) gives
(\ref{Tsp}).

Now we want to prove that $T^{\left(  s\right)  }$ can be represented by a
signed measure. We claim that
\begin{equation}
\int_{\rho}^{R}r^{-k}d\sigma_{k,l}^{\left(  s\right)  }\left(  r\right)
\leq\int_{\rho}^{R}r^{-k}d\mu_{k,l}\left(  r\right)  <\infty
\label{ChebyshevInequality}%
\end{equation}
for all $k\in\mathbb{N}_{0},$ $l=1,...,a_{k}.$ If $\sigma_{k,l}^{\left(
s\right)  }=\mu_{k,l}$ there is nothing to prove, so we can assume that the
cardinality of the support of $\mu_{k,l}$ is bigger than $s.$ After taking the
image measures under the map $\psi$ we see that we have to prove
\begin{equation}
\int_{\rho^{2}}^{R^{2}}t^{-\frac{1}{2}k}d\nu_{k,l}^{\left(  s\right)  }\left(
t\right)  \leq\int_{\rho^{2}}^{R^{2}}t^{-\frac{1}{2}k}d\mu_{k,l}^{\psi}\left(
t\right)  . \label{ChebyshevInequality2}%
\end{equation}
For $\rho>0$ this follows from the Chebyshev extremal property of the
Gau\ss--Jacobi measure (see e.g. \cite[Chapter 4, Theorem 1.1]{KrNu77})
applied to the function $f\left(  t\right)  :=t^{-\frac{1}{2}k}.$ The same
result works in the case $\rho=0$ but due to the singularity of $f$ we have to
use essentially the fact that all points of $\nu_{k,l}^{\left(  s\right)  }$
are in the open interval $\left(  0,R^{2}\right)  $ and to apply Remark $1.2$
in Chapter $4$ of \cite{KrNu77}.\footnote{It is curious that Stieltjes proved
that the Gau\ss--Jacobi quadrature measure solves a three--dimensional
spherically symmetric extremal problem with a singular function $f\left(
t\right)  =\frac{1}{\sqrt{t^{3}}},$ see the complete description in
\cite[Chapter 4.2, formula (2.6) ]{KrNu77}.} We have only to check that the
assumptions on $f$ are satisfied: $f$ is non-negative, $2s$ times
differentiable on the open interval $\left(  \rho^{2},R^{2}\right)  $, and for
$t\in\left(  \rho^{2},R^{2}\right)  $
\begin{equation}
\frac{d^{2s}f\left(  t\right)  }{dt^{2s}}=\left(  -1\right)  ^{2s}\frac{1}%
{2}k\left(  \frac{1}{2}k+1\right)  ....\left(  \frac{1}{2}k+2s-1\right)
t^{-\frac{1}{2}k-2s+1}>0. \label{eqFpositive}%
\end{equation}

By our assumption (\ref{condmu}) and by (\ref{ChebyshevInequality}) we can
apply Theorem \ref{ThmRepG} and (i) is proved. Property (ii) follows from
(\ref{Tsp}) which we have proved above.

Let us prove the uniqueness of $\sigma^{\left(  s\right)  }.$ Assume that
$\tau$ is a signed pseudo-positive measure with compact support, and with
properties (i) and (ii). Since $\tau$ is pseudo-positive there exists by
Proposition \ref{pseudopos} univariate measures $\tau_{k,l}$ such that
\begin{equation}
\int h\left(  t\right)  d\tau_{k,l}=\int h(\left|  x\right|  )Y_{k,l}\left(
x\right)  d\tau\left(  x\right)  \label{eqid1}%
\end{equation}
for any polynomially bounded continuous function $h$. Since $\Delta
^{2s}\left(  \left|  x\right|  ^{2j}Y_{k,l}\left(  x\right)  \right)  =0$ for
$j=0,...,2s-2,$ we infer that
\[
\int\left|  x\right|  ^{2j}Y_{k,l}\left(  x\right)  d\tau\left(  x\right)
=\int\left|  x\right|  ^{2j}Y_{k,l}\left(  x\right)  d\mu\left(  x\right)  .
\]
Hence $\int t^{2j}d\tau_{k,l}=\int t^{2j}d\mu_{k,l}$ for $j=0,...,2s-2,$ so
\begin{equation}
\int t^{j}d\tau_{k,l}^{\psi}=\int t^{j}d\mu_{k,l}^{\psi}. \label{identmark}%
\end{equation}
By property (i) the support $\tau_{k,l}$ has cardinality $\leq s,$ hence
$\tau_{k,l}^{\psi}$ has cardinality $\leq s.$ The uniqueness of the
Gau\ss-Jacobi quadrature shows that $\tau_{k,l}^{\psi}$ is equal to $\nu
_{k,l}^{\left(  s\right)  }$ which means that $\tau_{k,l}=\sigma_{k,l}.$ If
the support of $\mu_{k,l}$ has less than $s$ points then $\sigma
_{k,l}^{\left(  s\right)  }$ in our construction is defined to be $\mu_{k,l}.$
From (\ref{identmark}) one can derive that $\tau_{k,l}^{\psi}$ has the same
support set as $\mu_{k,l}^{\psi}$ and finally that $\tau_{k,l}=\mu
_{k,l}=\sigma_{k,l}^{\left(  s\right)  }$. Proposition \ref{soso} yields
$\tau=\sigma^{\left(  s\right)  }.$
\end{proof}

\begin{definition}
\label{DGJmeasure} The measure $\sigma^{\left(  s\right)  }$ constructed in
Theorem \ref{ThmMain} will be called the \textbf{polyharmonic Gau\ss-Jacobi
measure of order} $s$ for the measure $\mu.$
\end{definition}

The following is an analog to the theorem of Stieltjes about the convergence
of the univariate Gau\ss--Jacobi quadrature formulas.

\begin{theorem}
\label{TStieltjes}Let $0<R<\infty$ and let $\sigma^{\left(  s\right)  }$ be
the polyharmonic Gau\ss-Jacobi measure of order $s$ for the measure $\mu,$
obtained in Theorem \ref{ThmMain}. Then
\[
\int f\left(  x\right)  d\sigma^{\left(  s\right)  }\rightarrow\int f\left(
x\right)  d\mu\qquad\text{ for }s\rightarrow\infty
\]
holds for every function $f\in C\left(  B_{R}\right)  $.
\end{theorem}

\begin{proof}
Item (ii) of Theorem \ref{ThmMain} implies that for any polynomial $P$ the
convergence $T^{\left(  s\right)  }\left(  P\right)  \rightarrow P$ holds for
$s\longrightarrow\infty.$ Theorem 14.4.4 in \cite{Davis} shows that the
convergence $T^{\left(  s\right)  }\left(  f\right)  \rightarrow f$ carries
over to all continuous functions $f:B_{R}\rightarrow\mathbb{C}$ provided there
exists a constant $C>0$ such that for all natural numbers $s$ and all $f\in
C\left(  B_{R}\right)  $
\[
\left|  T^{\left(  s\right)  }\left(  f\right)  \right|  \leq C\max_{\left|
x\right|  \leq R}\left|  f\left(  x\right)  \right|  .
\]
But that is just estimate (\ref{eqkey23}) in the proof of Theorem
\ref{ThmRepG}. The proof is complete.
\end{proof}

Using the same techniques as in Theorem \ref{ThmMain} we may prove a
generalization of the Chebyshev extremal property of the Gau\ss--Jacobi quadrature:

\begin{theorem}
\label{Cextremal}Let $0\leq\rho<R<\infty$ and let $\mu$ be a pseudo-positive
signed measure with support in $A_{\rho,R}$ satisfying the summability
condition (\ref{condmu}) and let $\sigma^{\left(  s\right)  }$ be the
polyharmonic Gau\ss-Jacobi measure of order $s$. Let $f\in C^{2s}\left(
\mathbb{R}^{d}\right)  $ be such that for all $t\in\left(  \rho^{2}%
,R^{2}\right)  $ holds
\[
\frac{d^{2s}}{dt^{2s}}\left[  f_{k,l}\left(  \sqrt{t}\right)  t^{-\frac{1}%
{2}k}\right]  \geq0,
\]
for all $k\in\mathbb{N}_{0},$ $l=1,2,...,a_{k}$. Then the following
inequality
\[
\int f\left(  x\right)  d\sigma^{\left(  s\right)  }\leq\int f\left(
x\right)  d\mu
\]
holds.
\end{theorem}

Although the measures $\sigma_{k,l}^{\left(  s\right)  }$ are based on point
evaluations, it is clear that our approximation measures $\sigma^{\left(
s\right)  }$ are not point evaluations. The following gives a description of
the support of the polyharmonic Gau\ss--Jacobi measure $\sigma^{\left(
s\right)  }$ when only finitely many measures $\mu_{k,l}$ are non-zero:

\begin{proposition}
\label{Pfiniteonly} Suppose that $\sigma_{k,l},$ $k=1,...,k_{0},$
$l=1,...,a_{k},$ are non-negative measures with finite support in the open
interval $\left(  0,\infty\right)  $ and suppose that $\sigma_{k,l}=0$ for all
$k>k_{0}.$ Then the support of the representing measure of the functional
$T:\mathbb{C}\left[  x_{1},x_{2},...,x_{d}\right]  \rightarrow\mathbb{C}$
defined by (\ref{defTTT}) is contained in the union of finitely many spheres
with positive radius and with center $0.$
\end{proposition}

\begin{proof}
Let $S_{k,l}$ be the finite support of $\sigma_{k,l}$ and let $S$ be the union
of all sets $S_{k,l}$ with $k=1,...,k_{0},$ $l=1,...,a_{k}.$ Let
$\widetilde{S}=$ $\left\{  x\in\mathbb{R}^{d}:\left|  x\right|  \in S\right\}
.$ We show that the support of $\sigma$ is contained in $\widetilde{S}$.
Indeed, one can estimate as in (\ref{eqkey23})
\[
\left|  T\left(  P\right)  \right|  \leq\max_{x\in\widetilde{S}}\left|
P\left(  x\right)  \right|  \sqrt{\omega_{d-1}}\sum_{k=0}^{k_{0}}\sum
_{l=1}^{a_{k}}\int_{0}^{R}r^{-k}d\sigma_{k,l}\left(  t\right)  .
\]
By the Riesz representation theorem, $T$ can be represented by a signed
measure $\sigma$ which has support in the compact set $\widetilde{S}.$
\end{proof}

\begin{remark}
\label{Rminimalsupport}An interesting characteristic feature of the classical
Gau\ss--Jacobi quadrature measure is the minimality of its support among all
non--negative measures which are exact of the same degree $2s-1.$ One might
see above some analogy with this phenomenon if one considers $\mu_{k,l}\left(
r\right)  $ and $\sigma_{k,l}^{\left(  s\right)  }\left(  r\right)  $ as
measures defined on the space $K=\left\{  \left(  k,l,r\right)  :k\in
\mathbb{N}_{0},\ l=1,2,...,a_{k},\ r\in\left[  0,\infty\right)  \right\}  .$
\end{remark}

\section{Markov type error estimates\label{Serror}}

In this section we want to give an error estimate for our cubature formula.
The proof is based on the Markov estimate for the Gau\ss-Jacobi measure
provided in (\ref{mark}).

Let $s\in\mathbb{N}\cup\left\{  \infty\right\}  .$ For an open subset $U$ of
$\mathbb{R}^{d}$ we denote by $C^{s}\left(  U\right)  $ the space of all $f\in
C\left(  U\right)  $ which are continuously differentiable in $U$ up to the
order $s.$

\begin{theorem}
\label{PropError}Let $0\leq\rho<R<\infty$ and let $\psi:\left[  0,\infty
\right)  \rightarrow\left[  0,\infty\right)  $ be defined by $\psi\left(
t\right)  =t^{2}.$ Let $\mu$ be a pseudo-positive signed measure with support
in $A_{\rho,R}$ satisfying the summability condition (\ref{condmu}), and let
$\sigma^{\left(  s\right)  }$ be the polyharmonic Gau\ss-Jacobi measure of
order $s$. Define for every $f\in C\left(  A_{\rho,R}\right)  $ the error
functional
\[
E_{s}\left(  f\right)  :=\int f\left(  x\right)  d\mu\left(  x\right)  -\int
f\left(  x\right)  d\sigma^{\left(  s\right)  }\left(  x\right)  .
\]
If $f\in C^{2s}\left(  A_{\rho,R}^{o}\right)  \cap C\left(  A_{\rho,R}\right)
$ then $E_{s}\left(  f\right)  $ is lower equal to
\[
\frac{1}{\left(  2s\right)  !}\sum_{k=0}^{\infty}\sum_{l=1}^{a_{k}}\sup
_{\rho^{2}<\xi<R^{2}}\left|  \frac{d^{2s}}{dt^{2s}}\left[  f_{k,l}\left(
\sqrt{t}\right)  t^{-\frac{1}{2}k}\right]  \left(  \xi\right)  \right|
\int_{\rho^{2}}^{R^{2}}\left|  Q_{k,l}^{s}\left(  t\right)  \right|  ^{2}%
d\mu_{kl}^{\psi}.
\]
Here $Q_{k,l}^{s}\left(  t\right)  $ is the orthogonal polynomial of degree
$s$ with respect to the measure $\mu_{kl}^{\psi},$ normalized so that the
leading coefficient is equal to $1$; if the support of $\mu_{k,l}$ has less
than $s$ points, $Q_{k,l}^{s}$ is defined to be $0.$
\end{theorem}

\begin{proof}
Since $f\in C^{2s}\left(  A_{\rho,R}^{o}\right)  \cap C\left(  A_{\rho
,R}\right)  $ it is easy to see that the Laplace-Fourier coefficients
$f_{k,l}\in C^{2s}\left(  \rho,R\right)  \cap C\left[  \rho,R\right]  $. Let
$\mu_{k,l}$ and $\sigma_{k,l},$ $k\in\mathbb{N}_{0},$ $l=1,...,a_{k},$ and
$\sigma^{\left(  s\right)  }$ be as in Theorem \ref{ThmMain}. From the
definitions it follows
\[
E_{s}\left(  f\right)  =\sum_{k=0}^{\infty}\sum_{l=1}^{a_{k}}\int_{\rho}%
^{R}f_{k,l}\left(  r\right)  r^{-k}d\mu_{k,l}-\int_{\rho}^{R}f_{k,l}\left(
r\right)  r^{-k}d\sigma_{k,l}^{\left(  s\right)  }.
\]
Further $f_{k,l}\left(  r\right)  r^{-k}$ is integrable with respect to
$\mu_{k,l}$ since $f_{k,l}$ is continuous on $\left[  \rho,R\right]  $ and
condition (\ref{condmu}) holds. Let us fix the pair of indices $\left(
k,l\right)  .$ If the support of $\mu_{k,l}$ has less than $s$ points we know
that $\mu_{k,l}=\sigma_{k,l}^{\left(  s\right)  }$. So assume that the support
of $\mu_{k,l}$ has at least $s$ points. Then the support of $\mu_{k,l}^{\psi}
$ has at least $s$ points and in our construction $\nu_{k,l}^{\left(
s\right)  }$ is the Gau\ss-Jacobi measure of $\mu_{k,l}^{\psi}.$ Consequently
\begin{align*}
e\left(  f_{k,l}\right)   &  :=\int_{\rho}^{R}f_{k,l}\left(  r\right)
r^{-k}d\mu_{k,l}\left(  r\right)  -\int_{\rho}^{R}f_{k,l}\left(  r\right)
r^{-k}d\sigma_{k,l}^{\left(  s\right)  }\left(  r\right) \\
&  =\int_{\rho^{2}}^{R^{2}}f_{k,l}\left(  \sqrt{t}\right)  t^{-\frac{1}{2}%
k}d\mu_{k,l}^{\psi}\left(  t\right)  -\int_{\rho^{2}}^{R^{2}}f_{k,l}\left(
\sqrt{t}\right)  t^{-\frac{1}{2}k}d\nu_{k,l}^{\left(  s\right)  }\left(
t\right)  .
\end{align*}
By the proof of Markov's error estimate (\ref{mark}) given in (\cite{Davis})
one easily obtains with $g_{k,l}\left(  t\right)  :=f_{k,l}\left(  \sqrt
{t}\right)  t^{-\frac{1}{2}k}$ the inequality
\[
e\left(  f_{k,l}\right)  \leq\frac{1}{\left(  2s\right)  !}\sup_{\rho^{2}%
<\xi<R^{2}}\left|  g_{k,l}^{\left(  2s\right)  }\left(  \xi\right)  \right|
\int_{\rho^{2}}^{R^{2}}\left|  Q_{k,l}^{s}\left(  t\right)  \right|  ^{2}%
d\mu_{k,l}^{\psi}\left(  t\right)  .
\]
The proof is complete.
\end{proof}

In the following we want to give a Markov type error estimates for holomorphic
functions $f.$ We will need the following property which was observed in
\cite{Baouendi}:

\begin{lemma}
Let $f\in C^{\infty}\left(  B_{R}^{\circ}\right)  .$ Then $f_{k,l}\in
C^{\infty}\left[  0,R\right)  $ and $\frac{d^{m}}{dr^{m}}f_{k,l}\left(
0\right)  =0$ for $m=0,...,k-1.$
\end{lemma}

\begin{lemma}
Let $f$ be a holomorphic function on the open ball $B_{\tau}^{\mathbb{C}%
}:=\{w\in\mathbb{C}^{d}:\sum_{j=1}^{d}\left|  w_{j}\right|  ^{2}<\tau^{2}\}$
for $\tau>0.$ Let $f_{k,l}$ be the Laplace-Fourier coefficient of $f$ given by
(\ref{eqfourier}) and let $p_{k,l}\left(  t\right)  $ be defined by the
equation $f_{k,l}\left(  r\right)  =p_{k,l}\left(  r^{2}\right)  r^{k}$ for
$0<r<\tau.$ Let $\rho$ and $t$ satisfy $0\leq t^{2}<\rho<\tau.$ Then
\begin{equation}
\left|  \frac{d^{s}}{dt^{s}}p_{k,l}\left(  t\right)  \right|  \leq\sqrt
{\omega_{d}}\max_{u\in\left[  0,2\pi\right]  ,\theta\in\mathbb{S}^{d-1}%
}\left|  f\left(  e^{iu}\rho\theta\right)  \right|  \frac{\rho^{2-k}%
s!}{\left(  \rho^{2}-t\right)  ^{s+1}} \label{eqableitung}%
\end{equation}
hold for all $s=0,1,2,...$
\end{lemma}

\begin{proof}
Let $\theta\in\mathbb{S}^{d-1}.$ The map $\varphi_{\theta}:\left\{
z\in\mathbb{C}:\left|  z\right|  <\tau\right\}  \rightarrow B_{\tau
}^{\mathbb{C}}\left(  0\right)  $ defined by $\varphi_{\theta}\left(
z\right)  =z\theta$ is clearly holomorphic. Hence $f_{\theta}$ defined by
$f_{\theta}\left(  z\right)  =f\left(  z\theta\right)  =f\circ\varphi_{\theta
}\left(  z\right)  $ is holomorphic. It follows that $f_{k,l}\left(  z\right)
$ defined by
\begin{equation}
f_{k,l}\left(  z\right)  =\int_{\mathbb{S}^{d-1}}f\left(  z\theta\right)
Y_{k,l}\left(  \theta\right)  d\theta\label{deffkl}%
\end{equation}
is a holomorphic extension of $f_{k,l}$ to $\left\{  z\in\mathbb{C}:\left|
z\right|  <\tau\right\}  $. Cauchy's inequality shows that for $\rho=\left|
z\right|  $
\begin{align*}
\left|  f_{k,l}\left(  z\right)  \right|  ^{2}  &  \leq\int_{\mathbb{S}^{d-1}%
}\left|  f\left(  z\theta\right)  \right|  ^{2}d\theta\cdot\int_{\mathbb{S}%
^{d-1}}\left|  Y_{k,l}\left(  \theta\right)  \right|  ^{2}d\theta\\
&  \leq\omega_{d}\max_{u\in\left[  0,2\pi\right]  ,\theta\in\mathbb{S}^{d-1}%
}\left|  f\left(  e^{iu}\rho\theta\right)  \right|  ^{2}.
\end{align*}
The Cauchy estimates $\left|  g^{\left(  k\right)  }\left(  0\right)  \right|
\leq\frac{k!}{\rho^{k}}\max_{\left|  z\right|  =\rho}\left|  g\left(
z\right)  \right|  $ for a holomorphic function $g$ and the last estimate
imply for $0<\rho<\tau$%
\begin{equation}
\left|  \frac{d^{m+k}}{dz^{m+k}}f_{k,l}\left(  0\right)  \right|  \leq
\sqrt{\omega_{d}}\max_{u\in\left[  0,2\pi\right]  ,\theta\in\mathbb{S}^{d-1}%
}\left|  f\left(  e^{iu}\rho\theta\right)  \right|  \cdot\frac{\left(
k+m\right)  !}{\rho^{m+k}} \label{eqlast1}%
\end{equation}
Let us write $f_{k,l}\left(  z\right)  =\sum_{m=k}^{\infty}\frac{1}{m!}%
\frac{d^{m}}{dr^{m}}f_{k,l}\left(  0\right)  \cdot z^{m}$ for $\left|
z\right|  <\tau. $ It is known that $r^{-k}f_{k,l}\left(  r\right)  $ is an
even function, hence we can write
\[
p_{k,l}\left(  r^{2}\right)  =r^{-k}f_{k,l}\left(  r\right)  =\sum
_{m=0}^{\infty}\frac{1}{\left(  k+2m\right)  !}\frac{d^{2m+k}}{dr^{2m+k}%
}f_{k,l}\left(  0\right)  \cdot r^{2m}.
\]
Then for $t=r^{2}$ we have
\[
\frac{d^{s}}{dt^{s}}p_{k,l}\left(  t\right)  =\sum_{m=s}^{\infty}\frac
{1}{\left(  k+2m\right)  !}\frac{m!}{\left(  m-s\right)  !}\frac{d^{2m+k}%
}{dr^{2m+k}}f_{k,l}\left(  0\right)  \cdot t^{\left(  m-s\right)  }%
\]
Now (\ref{eqlast1}) implies
\[
\left|  \frac{d^{s}}{dt^{s}}p_{k,l}\left(  t\right)  \right|  \leq\sqrt
{\omega_{d}}\max_{u\in\left[  0,2\pi\right]  ,\theta\in\mathbb{S}^{d-1}%
}\left|  f\left(  e^{iu}\rho\theta\right)  \right|  \frac{1}{\rho^{k+2s}}%
\sum_{m=s}^{\infty}\frac{m!}{\left(  m-s\right)  !}\left(  \frac{t}{\rho^{2}%
}\right)  ^{m-s}.
\]
Since for $\left|  t\right|  <1$ we have $\sum_{m=s}^{\infty}\frac{m!}{\left(
m-s\right)  !}t^{m-s}=\frac{d^{s}}{dt^{s}}\sum_{m=0}^{\infty}t^{m}=\frac
{d^{s}}{dt^{s}}\frac{1}{1-t}=s!\left(  1-t\right)  ^{-s-1}$ a straightforward
computation yields (\ref{eqableitung}).
\end{proof}

Now combining Theorem \ref{PropError} and the last Lemma, we obtain the final estimate.

\begin{theorem}
Let $0<R<\infty$ and let $\mu$ be a pseudo-positive signed measure with
support in $B_{R}$ satisfying the summability condition (\ref{condmu}) and let
$\sigma^{\left(  s\right)  }$ be the polyharmonic Gau\ss-Jacobi measure of
order $s$. Then the following error estimate
\[
E\left(  f\right)  \leq\frac{\sqrt{\omega_{d}}\rho^{2}}{\left(  \rho^{2}%
-R^{2}\right)  ^{2s+1}}\max_{w\in\mathbb{C}^{n},\left|  w\right|  \leq\rho
}\left|  f\left(  w\right)  \right|  \sum_{k=0}^{\infty}\sum_{l=1}^{a_{k}%
}\frac{1}{\rho^{k}}\int_{0}^{R^{2}}\left|  Q_{k,l}^{s}\left(  t\right)
\right|  ^{2}d\mu_{k,l}^{\psi}\left(  t\right)
\]
holds for all functions $f:B_{R}\rightarrow\mathbb{C}$ which possess a
holomorphic extension to the complex ball $B_{\tau}^{\mathbb{C}}$ for $\tau>R
$ and for any $\rho$ with $R<\rho<\tau.$
\end{theorem}

\section{Polyharmonic Gau\ss--Jacobi cubature in the annulus \label{Sannulus}}

Let us imagine a function $f$ which is say holomorphic in the ball $B_{R}$
with some singularities in the smaller ball $B_{\rho}^{\circ}$, and one needs
to find the integral of the function on the annulus $A_{\rho,R}$ with an
estimate for the error of approximation. This example is a motivation to
consider polyharmonic Gau\ss--Jacobi cubatures in the annulus $A_{\rho,R}$
which generalize the construction of Section \ref{Scubature}.

By the results in Section \ref{Scubature} the polyharmonic Gau\ss-Jacobi
measure $\sigma^{\left(  s\right)  }$ of order $s$ has support in $A_{\rho,R}$
and it is exact on the space of all \emph{polynomials} $f$ such that
$\Delta^{2s}f=0.$ In the present section we will seek polyharmonic
Gau\ss--Jacobi cubatures which are exact on the larger space
\[
PH^{2s}\left(  A_{\rho,R}\right)  =\left\{  f\in C\left(  A_{\rho,R}\right)
\cap C^{2s}\left(  A_{\rho,R}^{\circ}\right)  :\Delta^{2s}f\left(  x\right)
=0\text{ for all }x\in A_{\rho,R}^{\circ}\right\}
\]
where $A_{\rho,R}^{\circ}$ is the interior of $A_{\rho,R}.$ Exactness with
respect to $PH^{2s}\left(  A_{\rho,R}\right)  $ is related to the expectation
that the integrals of functions with singularities in the inner open ball
$\left\{  \left|  x\right|  <\rho\right\}  $ will be better approximated.

It turns out that the problem can be solved in a way very similar to Theorem
\ref{ThmMain}. The proof is so far based on the ''generalized Gau\ss--Jacobi
quadratures for Chebyshev systems'' which have been developed mainly by A.
Markov, \cite[Chapter 4]{KrNu77}. We restrict our discussion to the case of
compact annulus which is less technical. Let us now formulate the result precisely:

\begin{theorem}
\label{TAnnulus}Let $0<\rho<R<\infty.$ Let $\mu$ be a pseudo-positive signed
measure with support in $A_{\rho,R}$ such that
\begin{equation}
\sum_{k=0}^{\infty}\sum_{l=1}^{a_{k}}\int_{\rho}^{R}r^{-k}d\mu_{k,l}\left(
r\right)  <\infty. \label{eqneu4}%
\end{equation}
Let $s$ be a natural number. Then there exists \textbf{unique pseudo-positive}
measure $\tau^{\left(  2s\right)  }$ with support in $A_{\rho,R}$ such that

(i) The cardinality of the support of each component measure $\tau
_{k,l}^{\left(  2s\right)  },$ $k\in\mathbb{N}_{0},$ $l=1,2,...,a_{k},$ is
$\leq2s.$

(ii) $\int f\left(  x\right)  d\tau^{\left(  2s\right)  }\left(  x\right)
=\int f\left(  x\right)  d\mu\left(  x\right)  $ for all $f\in PH^{2s}\left(
A_{\rho,R}\right)  .$
\end{theorem}

Let us compare the result with Theorem \ref{ThmMain}. In the latter case we
obtained a measure $\sigma^{\left(  s\right)  }$ with support in $A_{\rho,R}$
which is exact for all \emph{polynomials} $P$ with $\Delta^{2s}P=0.$ The
support of the component measure $\sigma_{k,l}^{\left(  s\right)  }$ has at
most $s$ points. In contrast, the component measure $\tau_{k,l}^{\left(
2s\right)  }$ of the solution $\tau^{\left(  2s\right)  }$ has a support of
cardinality $\leq2s$ which is twice bigger. This is caused by the fact that
$\tau^{\left(  2s\right)  }$ is exact on the larger subspace $PH^{2s}\left(
A_{\rho,R}\right)  .$

For the proof we will first need the representation of a polyharmonic function
in the annulus which is somewhat more sophisticated than that of a polynomial
as given in formula (\ref{gauss}). Let us introduce the operators
\begin{equation}
L_{\left(  k\right)  }:=\frac{d^{2}}{dr^{2}}+\frac{n-1}{r}\frac{d}{dr}%
-\frac{k\left(  n+k-2\right)  }{r^{2}} \label{Lk}%
\end{equation}
which may be written as
\begin{equation}
L_{\left(  k\right)  }f\left(  r\right)  =\frac{1}{r^{n+k-1}}\frac{d}%
{dr}\left[  r^{n+2k-1}\frac{d}{dr}\left[  \frac{1}{r^{k}}f\left(  r\right)
\right]  \right]  , \label{Lknicerepres}%
\end{equation}
see e.g. (10.18) in \cite{Koun00}. Let $L_{\left(  k\right)  }^{2s}$denote the
$2s$-th iterate of $L_{\left(  k\right)  }.$ These operators are the radial
part of the polyharmonic operators $\Delta^{2s}.$ A set of the $4s$ linearly
independent solutions of the equation
\begin{equation}
L_{\left(  k\right)  }^{2s}f\left(  r\right)  =0\qquad\text{for }r>0,
\label{Lkp}%
\end{equation}
can be explicitly constructed; e.g. for $k\geq4s$ a set of solutions is given
by
\begin{equation}
r^{2j-k}\qquad j=0,1,...,s-1;\quad r^{2j+k},\qquad j=0,1,...,2s-1.
\label{Lkp2}%
\end{equation}
For $k<4s$ one has to be careful with multiplicities, and we refer to
\cite{Koun00} for an explicit description.

We have the following result (see e.g. Theorem 10.39 in \cite{Koun00}):

\begin{proposition}
\label{Ppolyhrepresent}Let $h\in PH^{2s}\left(  A_{\rho,R}\right)  .$ Then the
Laplace--Fourier series
\begin{equation}
h\left(  x\right)  =\sum_{k=0}^{\infty}\sum_{l=1}^{d_{k}}h_{k,l}\left(
r\right)  Y_{k,l}\left(  \theta\right)  \label{eq70}%
\end{equation}
converges absolutely and uniformly on compact subsets of $A_{\rho,R}^{\circ}.
$ The Laplace--Fourier coefficients $h_{k,l}$ are solutions of (\ref{Lkp}).
\end{proposition}

In the following we will mimic the proof of Theorem \ref{ThmMain}. The
measures $\sigma_{k,l}^{\left(  s\right)  }$ in the proof of Theorem
\ref{ThmMain} had the feature that they were exact on the solutions
$1,r^{2},...,r^{4s-2}$. Proposition \ref{Ppolyhrepresent} shows that we need
now quadratures which are exact on the solutions of (\ref{Lkp}). This
motivates to recall the theory of A. Markov and M. Krein on quadratures for
Chebyshev systems.

\begin{definition}
\label{DChebyshevsystem} Let $u_{0},...,u_{N}$ be continuous functions on
$\left[  a,b\right]  .$ We say that they form a Chebyshev system of order
$N+1$ on $\left[  a,b\right]  $ if every non--trivial linear combination
$\ \sum_{j=0}^{N}\gamma_{j}u_{j}\left(  t\right)  $ has at most $N$ zeros on
$\left[  a,b\right]  ,$ i.e. the determinant
\begin{equation}
\det\left(  u_{j}\left(  t_{i}\right)  \right)  _{i,j=0}^{N} \label{det}%
\end{equation}
is not zero on $\left[  a,b\right]  .$ The system $u_{0},...,u_{N}$ is $T_{+}$
on $\left[  a,b\right]  $ if $\det\left(  u_{j}\left(  t_{i}\right)  \right)
_{i,j=0}^{N}>0$ holds for every choice of $t_{j}\in$ $\left[  a,b\right]  $
with $t_{0}<t_{1}<...<t_{N}$.
\end{definition}

Note that the definition of a $T_{+}$-system depends on the order of the
functions $u_{0},...,u_{N}$.

The following theorem is the generalization of the Gau\ss-Jacobi quadratures
for Chebyshev systems (see Theorem $1.1$ and Remark $1.2$ in Chapter $4,$ and
Theorem $1.4$ in Chapter $2$ in \cite{KrNu77}).

\begin{theorem}
\label{TMarkov}Let $N=2s-1$ for a natural number $s$ and let $\sigma$ be a
non--negative measure on $\left[  a,b\right]  $ with cardinality of the
support $>s.$ Let the continuous functions $u_{0},...,u_{N}$ be a Chebyshev
system on the interval $\left[  a,b\right]  $ and assume that $u_{N+1}%
:=\Omega$ is a continuous function on $\left[  a,b\right]  .$ If
$u_{0},...,u_{N+1}$ is a $T_{+}$-system on $\left[  a,b\right]  $ then there
exists a unique measure $\sigma^{\left(  s\right)  }$ with support of
cardinality $s$ such that
\begin{equation}
\int_{a}^{b}u_{j}\left(  t\right)  d\sigma\left(  t\right)  =\int_{a}^{b}%
u_{j}\left(  t\right)  d\sigma^{\left(  s\right)  }\left(  t\right)
\qquad\text{for }j=0,1,...,N. \label{moments}%
\end{equation}
The support of $\sigma^{\left(  s\right)  }$ is contained in the open interval
$\left(  a,b\right)  $.
\end{theorem}

The measure $\sigma^{\left(  s\right)  }$ is called in \cite{KrNu77} the
''lower chief representation'' which is also very natural to be called
\emph{Gau\ss--Jacobi--Markov quadrature}, and we will use this name further.

A second major result in the Krein-Markov theory is the extremal property of
the truncated moment problem due to Chebyshev, Markov and Stieltjes.

\begin{theorem}
\label{TMarkov2}With the notations and assumptions of Theorem \ref{TMarkov}
let $\sigma^{\left(  s\right)  }$ be the Gau\ss--Jacobi--Markov quadrature of
$\sigma$. The measure $\sigma^{\left(  s\right)  }$ attains the minimum in the
problem
\begin{equation}
\min_{\nu}\int_{a}^{b}\Omega\left(  t\right)  d\nu\left(  t\right)
\label{extremalproblem}%
\end{equation}
where $\nu$ ranges over all non--negative measures $\nu$ such that
\[
\int_{a}^{b}u_{j}\left(  t\right)  d\nu\left(  t\right)  =\int_{a}^{b}%
u_{j}\left(  t\right)  d\sigma\left(  t\right)  \qquad\text{for }j=0,1,...,N.
\]
\end{theorem}

We return now to our case of polyharmonic cubatures on annuli. At first we note

\begin{proposition}
\label{PsystemR} Any linear independent system $R_{k,j}^{s},j=1,2,...,4s,$ of
solutions of (\ref{Lkp}) is a Chebyshev system of order $2s$ on every interval
$\left[  a,b\right]  $ with $a>0.$ For $k\geq2s$ the system of solutions in
(\ref{Lkp2}) is a $T_{+}-$system on $\left[  a,b\right]  $ with $a>0.$
\end{proposition}

The proof follows from the results in Section II.5 and Theorem II.5.2 in
\cite{KrNu77} and uses the representation (\ref{Lknicerepres}). In view of the
Krein-Markov theory we need the following stronger result:

\begin{proposition}
\label{Pk>=4s}Let $k\in\mathbb{N}_{0},$ $l=1,2,...,a_{k}$ be fixed, with
$k\geq4s$, and define $u_{4s}=\Omega\equiv1.$ If we denote by $u_{0}%
,...,u_{4s-1}$ the system of solutions in (\ref{Lkp2}) then $u_{0},...,u_{4s}$
is a $T_{+}-$system on $\left[  a,b\right]  $ for $a>0.$
\end{proposition}

\begin{proof}
By the example in \cite[Chapter II, Section 2.1/c)]{KrNu77} the system
\[
\left\{  e^{t\left(  2j-k\right)  },\quad j=0,1,...,2s-1;\quad e^{0};\quad
e^{t\left(  2j+k\right)  },\quad j=0,1,...,2s-1\right\}
\]
is a $T_{+}-$system since the numbers $2j-k$ and $2j+k$ are all different due
to $k\geq4s.$ Then the reordered system
\[
\left\{  e^{t\left(  2j-k\right)  },\quad j=0,1,...,2s-1;\quad e^{t\left(
2j+k\right)  },\quad j=0,1,...,2s-1;\quad e^{0}\right\}
\]
has the same determinant sign of (\ref{det}) as the above. By a change of the
variable $r=e^{t},$ one concludes that the system $u_{0},...,u_{4s}$ is a
$T_{+}-$system.
\end{proof}

Now we are prepared to make the proof.

\begin{proof}
[Proof of Theorem \ref{TAnnulus}]Fix a pair of indices $\left(  k,l\right)  $
with $k\in\mathbb{N}_{0},$ $l=1,2,...,a_{k}$, and let $\mu_{k,l}$ be the
component measure. If the support of $\mu_{k,l}$ has less than $2s$ points,
put $\tau_{k,l}^{\left(  2s\right)  }=\mu_{k,l}.$ Assume now it has more than
$2s$ points. Let $u_{0},...,u_{4s-1}$ be the system (\ref{Lkp2}). By
Markov--Krein's Theorem \ref{TMarkov} applied to $\sigma:=\mu_{k,l}$ there
exists a Gau\ss--Jacobi--Markov measure $\tau_{k,l}^{\left(  2s\right)  }$
with support in $\left(  \rho,R\right)  $, and its support has cardinality
$2s.$ An essential point is to prove that (at least) for sufficiently large
$k$ one has
\begin{equation}
\int_{\rho}^{R}1d\tau_{k,l}^{\left(  2s\right)  }\left(  r\right)  \leq
\int_{\rho}^{R}1d\mu_{k,l}\left(  r\right)  . \label{ChebyshevInequalityAnn}%
\end{equation}
For $k\geq4s$ this follows immediately from Markov--Krein's Theorem
\ref{TMarkov2} by means of Proposition \ref{Pk>=4s}.

Further we proceed as in the proof of Theorem \ref{ThmMain}. We want to define
a functional $T^{\left(  2s\right)  }$ on $C\left(  A_{\rho,R}\right)  $ by
putting
\begin{equation}
T^{\left(  2s\right)  }\left(  f\right)  :=\sum_{k=0}^{\infty}\sum
_{l=1}^{a_{k}}\int_{\rho}^{R}f_{k,l}\left(  r\right)  d\tau_{k,l}^{\left(
2s\right)  }\left(  r\right)  \label{neurepp}%
\end{equation}
for $f\in C\left(  A_{\rho,R}\right)  $ where $f_{k,l}$ are its
Laplace--Fourier coefficients. Indeed, by using the standard estimate for the
Laplace-Fourier coefficients the inequality
\[
\left|  T^{\left(  2s\right)  }\left(  f\right)  \right|  \leq C\max_{\rho
\leq\left|  x\right|  \leq R}\left|  f\left(  x\right)  \right|  \sum
_{k=0}^{\infty}\sum_{l=1}^{a_{k}}\int_{\rho}^{R}d\tau_{k,l}^{\left(
2s\right)  }\left(  r\right)
\]
is easily established for all $f\in C\left(  A_{\rho,R}\right)  $. Now with
(\ref{ChebyshevInequalityAnn}) and our assumption (\ref{eqneu4}) it follows
that $T^{\left(  2s\right)  }$ is well-defined. We may apply the Riesz
representation theorem and obtain a representing measure, denoted by
$\tau^{\left(  2s\right)  }$, with support in $A_{\rho,R}$. Since the constant
function is in $C\left(  A_{\rho,R}\right)  $ it is clear that $\tau^{\left(
2s\right)  }$ is a finite measure. Let us remark that due to our assumption
(\ref{eqneu4}) the identity
\[
\int_{A_{\rho,R}}f\left(  x\right)  d\mu=\sum_{k=0}^{\infty}\sum_{l=1}^{a_{k}%
}\int_{\rho}^{R}f_{k,l}\left(  r\right)  d\mu_{k,l}\left(  r\right)
\]
holds for all $f\in C\left(  A_{\rho,R}\right)  ,$ since the right hand side
defines a continuous functional on $C\left(  A_{\rho,R}\right)  $ which agrees
with $\mu$ on the dense subspace $C^{\times}\left(  A_{\rho,R}\right)  $. Due
to the exactness property of all measures $\tau_{k,l}^{\left(  2s\right)  }$
and the representation (\ref{neurepp}), it follows that $\tau^{\left(
2s\right)  }$ satisfies ii) of Theorem \ref{TAnnulus} for all $f\in
PH^{2s}\left(  A_{\rho,R}\right)  $.

The pseudo-positivity of $\tau^{\left(  2s\right)  }$ follows from Corollary
\ref{ppp} since $T^{\left(  2s\right)  }$ is clearly pseudo-positive definite.
The uniqueness of $\tau^{\left(  2s\right)  }$ follows from the uniqueness of
the Gau\ss--Jacobi--Markov measure $\tau_{k,l}^{\left(  2s\right)  }$ as in
the proof of Theorem \ref{ThmMain}.
\end{proof}

\section{Polyharmonic Gau\ss--Jacobi cubature in the cylinder (periodic strip)
\label{Sstrip}}

The concept of pseudo-positivity which we have studied so far depends on the
expansion of the polyharmonic functions in Laplace-Fourier series which uses
the rotational symmetry of the ball and the annulus. The polyharmonic
Gau\ss-Jacobi cubature in the ball was defined respectively by the application
of the Gau\ss-Jacobi quadrature to the Laplace-Fourier coefficients. It is
natural to extend this concept on expansions available in other domains with
symmetries, and we will do so for the case of the cylinder (which may be
considered also as a periodic strip), where the Fourier series is the natural expansion.

Let now $-\infty\leq a<b\leq\infty.$ We consider functions $f:\left[
a,b\right]  \times\mathbb{R}^{d-1}\rightarrow\mathbb{C}$ depending on
$x=\left(  t,y\right)  $ in the \emph{strip }$\left[  a,b\right]
\times\mathbb{R}^{d-1}$ which are $2\pi$-periodic\footnote{The case of
non--periodic functions on the strip $\left[  a,b\right]  \times
\mathbb{R}^{d-1}$ is similar but needs more care.} with respect to the
variable $y,$ i.e. satisfy the equality
\[
f\left(  t,y+2\pi k\right)  =f\left(  t,y\right)  \qquad\text{ for all }%
t\in\mathbb{R},\text{ }y\in\mathbb{R}^{d-1},\ k\in\mathbb{Z}^{d-1}.
\]
Let us introduce the cylinder
\[
X=\left[  a,b\right]  \times\mathbb{T}^{d-1}%
\]
where $\mathbb{T}^{d-1}=\mathbb{S}^{d-1}\ $is the $d-1$ dimensional torus. One
may interpret the space $X$ also as a \emph{periodic strip}. The space of
$2\pi-$periodic in $y$ functions coincides with the space $C\left(  X\right)
.$ By $PH^{s}\left(  \left[  a,b\right]  \times\mathbb{T}^{d-1}\right)  $ we
denote the space of functions $f:X\rightarrow\mathbb{C}$ which are
\emph{polyharmonic of order} $s$ on $X,$ i.e. the space of functions $f\in X$
such that the corresponding $2\pi-$periodic in $y$ function $f\in C\left(
\left[  a,b\right]  \times\mathbb{R}^{d-1}\right)  $ is polyharmonic of order
$s $ in $\left[  a,b\right]  \times\mathbb{R}^{d-1}$ of the variables $x.$

If we fix $t\in\left[  a,b\right]  ,$ the equality
\begin{equation}
f\left(  t,y\right)  =\sum_{k\in\mathbb{Z}^{d-1}}f_{k}\left(  t\right)
e^{i\left\langle k,y\right\rangle } \label{fty}%
\end{equation}
is the Fourier series expansion of $f\in C\left(  X\right)  $ where the
Fourier coefficients $f_{k}\left(  t\right)  ,$ $k\in\mathbb{Z}^{d-1},$ of
$f\left(  t,\cdot\right)  $ are given by
\begin{equation}
f_{k}\left(  t\right)  =\frac{1}{\left(  2\pi\right)  ^{d-1}}\int
_{\mathbb{T}^{d-1}}f\left(  t,y\right)  e^{-i\left\langle k,y\right\rangle
}dy=\frac{1}{\left(  2\pi\right)  ^{d-1}}\int_{0}^{2\pi}....\int_{0}^{2\pi
}f\left(  t,y\right)  e^{-i\left\langle k,y\right\rangle }dy. \label{fk}%
\end{equation}

In the following we want to construct polyharmonic Gau\ss-Jacobi cubatures for
measures $\mu$ defined on $\mathbb{R}\times\mathbb{T}^{d-1},$ or equivalently,
defined on $\mathbb{R}^{d}$ and $2\pi-$periodic with respect to the variable
$y$.

Next we introduce pseudo-positivity in this setting:

\begin{definition}
\label{Dstrippseudopos} A measure $\mu$ on $\mathbb{R}\times\mathbb{T}^{d-1}$
with support in the cylinder $X$ is \textbf{pseudo-positive }(in order to
avoid mixing with the notion of pseudo--positivity introduced in Section
\ref{Smomentproblem} we will say sometimes \textbf{pseudo--positive on the
cylinder }$X$), if for each $k\in\mathbb{Z}^{d-1}$ and for each non-negative
continuous function $h:\mathbb{R}\rightarrow\left[  0,\infty\right)  $ with
compact support the inequality
\[
\int_{X}h\left(  t\right)  e^{i\left\langle k,y\right\rangle }d\mu\left(
t,y\right)  \geq0
\]
holds.
\end{definition}

For every $k\in\mathbb{Z}^{d-1},$ one may apply the Riesz representation
theorem to prove the existence of a unique non-negative measure $\mu_{k}$ on
$\mathbb{R}$ such that
\begin{equation}
\int_{-\infty}^{\infty}h\left(  t\right)  d\mu_{k}\left(  t\right)  =\int
_{X}h\left(  t\right)  e^{i\left\langle k,y\right\rangle }d\mu\left(
t,y\right)  \label{defneumu}%
\end{equation}
holds for any $h\in C_{c}\left(  \mathbb{R}\right)  .$ If we assume that $\mu$
has support in $X=\left[  a,b\right]  \times\mathbb{T}^{d-1}\mathbb{\ }$with
$-\infty<a<b<\infty$ then it is clear that (\ref{defneumu}) holds for all
$h\in C\left(  \mathbb{R}\right)  .$

Now we are going to prove the existence of polyharmonic Gau\ss--Jacobi
cubature for the case of the cylinder.

\begin{theorem}
\label{Tmainstrip}Let $-\infty<a<b<\infty.$ Let $\mu$ be a finite, signed
measure with support in the cylinder $X=\left[  a,b\right]  \times
\mathbb{T}^{d-1}\mathbb{\ }$which is pseudo-positive on $X$. Suppose that
\[
C:=\sum_{k\in\mathbb{Z}^{d-1}}\int_{X}e^{i\left\langle k,y\right\rangle }%
d\mu\left(  t,y\right)  <\infty.
\]
Then for each natural number $s$ there exists a unique finite signed measure
$\sigma^{\left(  2s\right)  }$ with support in $X$ such that

(i) The support of each component measure $\sigma_{k}^{\left(  2s\right)
},k\in\mathbb{Z}^{d-1}$ is in $\left[  a,b\right]  $ and has cardinality $\leq2s.$

(ii) $\int fd\mu=\int fd\sigma^{\left(  2s\right)  }$ for all functions $f\in
PH^{2s}\left(  X\right)  .$

(iii) $\sigma^{\left(  2s\right)  }$ is pseudo-positive on the cylinder $X$.
\end{theorem}

\begin{proof}
Let us give a characterization of the space $PH^{2s}\left(  X\right)  .$ If
$f\in PH^{2s}\left(  X\right)  $ then $f$ is representable in the Fourier
series (\ref{fty}) where for every $k\in\mathbb{Z}^{n-1}$ the function
$f_{k}\left(  t\right)  $ is a $C^{\infty}$-solution to the equation
\[
\left(  \frac{d^{2}}{dt^{2}}-k^{2}\right)  ^{2s}g\left(  t\right)  =0,
\]
cf. Theorem $9.3$ in \cite{Koun00}. All solutions of the latter equation are
linear combinations of the following functions
\begin{align*}
u_{j}\left(  t\right)   &  =t^{j}e^{-\left|  k\right|  t}\qquad\text{for
}j=0,1,...,2s-1,\\
u_{2s+j}\left(  t\right)   &  =t^{j}e^{\left|  k\right|  t}\qquad\text{for
}j=0,1,...,2s-1.
\end{align*}
Further define the function $u_{4s}\equiv1.$ Then the system of functions
$u_{0},...,u_{4s}$ is a $T_{+}$ system -- this follows from example c) in
Chapter 2.2 of \cite{KrNu77}.

Let $\mu_{k}$ be the non-negative measure defined by (\ref{defneumu}). If the
support of $\mu_{k}$ has less or equal than $2s$ points we define $\sigma
_{k}^{\left(  2s\right)  }:=\mu_{k}.$ If the support of $\mu_{k}$ has more
than $2s$ points, there exists according to Theorem \ref{TMarkov} a
non-negative measure $\sigma_{k}^{\left(  2s\right)  }$ such that
\[
\int_{a}^{b}u_{j}\left(  t\right)  d\sigma_{k}^{\left(  2s\right)  }\left(
t\right)  =\int_{a}^{b}u_{j}\left(  t\right)  d\mu_{k}\left(  t\right)
\qquad\text{for }j=0,...,4s-1,
\]
and the support of $\sigma_{k}^{\left(  2s\right)  }$ has $\leq2s$ points and
lies in $\left[  a,b\right]  .$ From the Krein--Markov Theorem \ref{TMarkov2}
it follows
\begin{equation}
\int_{\mathbb{R}}1d\sigma_{k}\left(  t\right)  \leq\int_{\mathbb{R}}1d\mu
_{k}\left(  t\right)  . \label{ChebyshevInequalityStrip}%
\end{equation}

Following the usual scheme, we want to define the functional $T^{\left(
2s\right)  }$ on\qquad$\allowbreak C\left(  X\right)  $ by putting
\[
T^{\left(  2s\right)  }\left(  f\right)  :=\frac{1}{\left(  2\pi\right)
^{n-1}}\sum_{k\in\mathbb{Z}^{n-1}}\int_{\mathbb{R}}f_{k}\left(  t\right)
d\sigma_{k}^{\left(  2s\right)  }\left(  t\right)  ,
\]
We have to show that the functional $T^{\left(  2s\right)  }$ is well-defined.
Note that
\begin{equation}
\left|  T^{\left(  2s\right)  }\left(  f\right)  \right|  \leq\frac{1}{\left(
2\pi\right)  ^{n-1}}\sum_{k\in\mathbb{Z}^{n-1}}\max_{t\in\left[  a,b\right]
}\left|  f_{k}\left(  t\right)  \right|  \int_{\mathbb{R}}1d\sigma_{k}\left(
t\right)  . \label{eqlonsum2}%
\end{equation}
Moreover it is clear that
\[
\max_{t\in\left[  a,b\right]  }\left|  f_{k}\left(  t\right)  \right|
\leq\max_{t\in\left[  a,b\right]  }\max_{y\in\mathbb{T}^{d-1}}\left|  f\left(
t,y\right)  \right|  .
\]
By (\ref{ChebyshevInequalityStrip}) we obtain
\[
\left|  T^{\left(  2s\right)  }\left(  f\right)  \right|  \leq\max
_{t\in\left[  a,b\right]  }\max_{y\in\mathbb{T}^{d-1}}\left|  f\left(
t,y\right)  \right|  \sum_{k\in\mathbb{Z}^{n-1}}\int_{\mathbb{R}}1d\mu
_{k}\left(  t\right)  ,
\]
Thus the functional $T^{\left(  2s\right)  }\left(  f\right)  $ is
well-defined on $C\left(  X\right)  ,$ and by the Riesz representation theorem
there exists a signed representing measure with support in $X$, denoted by
$\sigma^{\left(  2s\right)  }.$ Arguments similar to those presented at the
end of the proof of Theorem \ref{TAnnulus} show that $\sigma^{\left(
2s\right)  }$ is pseudo-positive, unique and that the exactness property (ii)
is satisfied. The details are omitted.
\end{proof}

\section{Examples and miscellaneous results\label{Sexamples}}

In this section we provide some examples and results on pseudo--positive
measures which throw more light on these new notions.

\subsection{The univariate case\label{Sunivariate}}

It is instructive to consider the univariate case of our theory: then $d=1,$
$\mathbb{S}^{0}=\left\{  -1,1\right\}  ,$ and the normalized measure is
$\omega_{0}\left(  \theta\right)  =\frac{1}{2}$ for all $\theta\in
\mathbb{S}^{0}.$ The harmonic polynomials are the linear functions, their
basis are the two functions defined by $Y_{0}\left(  x\right)  =1$ and
$Y_{1}\left(  x\right)  =x$ for all $x\in\mathbb{R}.$ The following is now
immediate from the definitions:

\begin{proposition}
Let $d=1.$ A functional $T:\mathbb{C}\left[  x\right]  \rightarrow\mathbb{C\ }
$ is pseudo-positive definite if and only if $T\left(  p^{\ast}\left(
x^{2}\right)  p\left(  x^{2}\right)  \right)  \geq0$ and $T\left(  xp^{\ast
}\left(  x^{2}\right)  p\left(  x^{2}\right)  \right)  \geq0$ for all
$p\in\mathbb{C}\left[  x\right]  .$
\end{proposition}

It follows from the last proposition that a Stieltjes moment sequence is
always pseudo-positive definite; by definition the functional $T:\mathbb{C}%
\left[  x\right]  \rightarrow\mathbb{C\ }$has the stronger property that
$T\left(  q^{\ast}\left(  x\right)  q\left(  x\right)  \right)  \geq0$ and
$T\left(  xq^{\ast}\left(  x\right)  q\left(  x\right)  \right)  \geq0$ for
all $q\in\mathbb{C}\left[  x\right]  .$ Below we give an example of a
pseudo-positive definite functional which is not positive definite, in
particular it does not define a Stieltjes moment sequence.

As pointed out in \cite[Chapter 4.1]{StWe71}, the Laplace--Fourier expansion
of $f$ is given by
\[
f\left(  r\theta\right)  =f_{0}\left(  r\right)  Y_{0}\left(  \theta\right)
+f_{1}\left(  r\right)  Y_{1}\left(  \theta\right)
\]
for $x=r\theta$ with $r=\left|  x\right|  $ and $\theta\in\mathbb{S}^{0},$
where
\begin{align*}
f_{0}\left(  r\right)   &  =\int_{\mathbb{S}^{0}}f\left(  r\theta\right)
Y_{0}\left(  \theta\right)  d\omega_{0}\left(  \theta\right)  =\frac{f\left(
r\right)  +f\left(  -r\right)  }{2}\\
f_{1}\left(  r\right)   &  =\int_{\mathbb{S}^{0}}f\left(  r\theta\right)
Y_{1}\left(  \theta\right)  d\omega_{0}\left(  \theta\right)  =\frac{f\left(
r\right)  -f\left(  -r\right)  }{2}%
\end{align*}
are the usual even and odd functions.

\begin{example}
Let $\sigma$ be a non-negative finite measure on the interval $\left[
a,b\right]  $ with $a>0.$ Then the functional $T:\mathbb{C}\left[  x\right]
\rightarrow\mathbb{C\ }$ defined by
\[
T\left(  f\right)  =\int_{a}^{b}f\left(  x\right)  d\sigma-\int_{a}%
^{b}f\left(  -x\right)  d\sigma
\]
is pseudo-positive definite but not positive definite.
\end{example}

\begin{proof}
Let $f_{1}$ be as above. Then $T\left(  f\right)  =2\int_{a}^{b}f_{1}\left(
r\right)  d\sigma,$ so $T$ is pseudo-positive definite by Proposition
\ref{PropTTT}. Since $T\left(  1\right)  =0$ and $T\neq0$ it is clear that $T$
is not positive definite.
\end{proof}

\subsection{A criterion for pseudo-positivity}

The following is a simple criterion for pseudo-positivity:

\begin{proposition}
\label{ThmLaplace}Let $\mu$ be a signed moment measure on $\mathbb{R}^{d}.$
Assume that $\mu$ has a density $w\left(  x\right)  $ with respect to the
Lebesgue measure $dx$ such that $\theta\longmapsto w\left(  r\theta\right)  $
is in $L^{2}\left(  \mathbb{S}^{d-1}\right)  $ for each $r>0.$ If the
Laplace-Fourier coefficients of $w,$
\[
w_{k,l}\left(  r\right)  :=\int_{\mathbb{S}^{d-1}}w\left(  r\theta\right)
Y_{k,l}\left(  \theta\right)  d\theta
\]
are non-negative then $\mu$ is pseudo-positive and
\begin{align}
d\mu_{k,l}\left(  r\right)   &  =r^{k+d-1}w_{k,l}\left(  r\right)
,\label{powermmuk}\\
\int_{0}^{\infty}r^{-k}d\mu_{k,l}\left(  r\right)   &  =\int_{0}^{\infty
}w_{k,l}\left(  r\right)  \cdot r^{d-1}dr \label{powermmuk2}%
\end{align}
if the last integral exists. The measures $\mu_{k,l}$ are defined by means of
equality (\ref{eqlim}).
\end{proposition}

\begin{proof}
Since $\mu$ has a density $w\left(  x\right)  $ we can use polar coordinates
to obtain for $f\in C_{pol}\left(  \mathbb{R}^{d}\right)  $
\begin{equation}
\int_{\mathbb{R}^{d}}fd\mu=\int_{\mathbb{R}^{d}}f\left(  x\right)  w\left(
x\right)  dx=\int_{0}^{\infty}\int_{\mathbb{S}^{d-1}}f\left(  r\theta\right)
w\left(  r\theta\right)  r^{d-1}d\theta dr. \label{eqPolar}%
\end{equation}
For any $h\in C_{pol}\left[  0,\infty\right)  $ we put $f\left(  x\right)
=h\left(  \left|  x\right|  \right)  Y_{k,l}\left(  x\right)  ,$ then we
obtain
\begin{equation}
\int_{\mathbb{R}^{d}}h\left(  \left|  x\right|  \right)  Y_{k,l}\left(
x\right)  d\mu=\int_{0}^{\infty}\int_{\mathbb{S}^{d-1}}h\left(  r\right)
r^{k+d-1}Y_{k,l}\left(  \theta\right)  w\left(  r\theta\right)  d\theta dr.
\label{eqPolar2}%
\end{equation}
Since $\theta\longmapsto w\left(  r\theta\right)  $ is in $L^{2}\left(
\mathbb{S}^{d-1}\right)  $, we know that $w_{k,l}\left(  r\right)
=\int_{\mathbb{S}^{d-1}}w\left(  r\theta\right)  Y_{k,l}\left(  \theta\right)
d\theta.$ Hence, by the definition of $\mu_{k,l},$ we obtain
\begin{equation}
\int_{0}^{\infty}h\left(  r\right)  d\mu_{k,l}:=\int_{\mathbb{R}^{d}}h\left(
\left|  x\right|  \right)  Y_{k,l}\left(  x\right)  d\mu=\int_{0}^{\infty
}h\left(  r\right)  w_{k,l}\left(  r\right)  r^{k+d-1}dr. \label{eqLaplacemu}%
\end{equation}
Thus the measure $\mu$ is pseudo-positive, and (\ref{powermmuk}) follows. Let
us prove (\ref{powermmuk2}): we define the cut--off functions $h_{m}\in
C_{pol}\left[  0,\infty\right)  $ such that $h_{m}\left(  t\right)  =t^{-k}$
for $t\geq1/m$ and such that $h_{m}\leq h_{m+1}.$ Now use (\ref{eqLaplacemu})
and the monotone convergence theorem to obtain (\ref{powermmuk2}).
\end{proof}

The next example addresses the question of whether there is a relationship
among the supports of the Gau\ss-Jacobi quadratures $\sigma_{k,l}$ in Theorem
\ref{ThmMain}:

\begin{proposition}
Let $s$ be a natural number. Then there exists a pseudo-positive measure $\mu$
with support in the unit ball such that the component measures $\sigma
_{k,1}^{\left(  s\right)  }$ of the polyharmonic Gau\ss-Jacobi cubature in
Theorem \ref{ThmMain} of order $s$ have identical supports.
\end{proposition}

\begin{proof}
Let $0<a<b$ and consider the density
\[
w\left(  re^{it}\right)  =\sum_{k=1}^{\infty}1_{\left[  a,b\right]  }\left(
r\right)  \frac{a^{k+d-1}}{r^{k+d-1}}\frac{1}{k^{2}}Y_{k,1}\left(
\theta\right)
\]
where $1_{\left[  a,b\right]  }$ is the indicator function of $\left[
a,b\right]  .$ Then $w$ induces a pseudo-positive measure and $\int
pd\mu_{k,1}=\frac{a^{k+d-1}}{k^{2}}\int_{a^{2}}^{b^{2}}p\left(  t\right)  dt$
according to (\ref{eqLaplacemu}). It follows that for all $k\in\mathbb{N}$ the
orthogonal polynomials of degree $s$ associated with $\mu_{k,1}$ are identical
up to a factor. Hence the supports of the measures $\sigma_{k,1}^{\left(
2s\right)  }$ are identical for all $k\in\mathbb{N}.$
\end{proof}

\subsection{The two--dimensional case\label{S2d}}

Let us consider the case $d=2,$ and take the usual orthonormal basis of solid
harmonics, defined by $Y_{0}\left(  e^{it}\right)  =\frac{1}{2\pi}$ and
\begin{equation}
Y_{k,1}\left(  re^{it}\right)  =\frac{1}{\sqrt{\pi}}r^{k}\cos kt\text{ and
}Y_{k,2}\left(  re^{it}\right)  =\frac{1}{\sqrt{\pi}}r^{k}\sin kt\text{ for
}k\in\mathbb{N}. \label{eqY}%
\end{equation}
We define a density $w^{\left(  \alpha\right)  }:\mathbb{R}^{n}\rightarrow
\left[  0,\infty\right)  $, depending on parameter $\alpha>0$, by
\begin{align*}
w^{\left(  \alpha\right)  }\left(  re^{it}\right)   &  :=\left(  1-r^{\alpha
}\right)  P\left(  re^{it}\right)  \qquad\text{ for }0\leq r<1\\
w^{\left(  \alpha\right)  }\left(  re^{it}\right)   &  =0\qquad\text{for
}r\geq1;
\end{align*}
here the function $P\left(  re^{it}\right)  $ is the Poisson kernel for $0\leq
r<1$ given by (see e.g. 5.1.16 in \cite[p. 243]{AAR99})
\begin{equation}
P\left(  re^{it}\right)  :=\frac{1-r^{2}}{1-2r\cos t+r^{2}}=1+\sum
_{k=1}^{\infty}2r^{k}\cos kt. \label{eqPoisson}%
\end{equation}
By Proposition \ref{ThmLaplace}, the measure $d\mu^{\alpha}:=w^{\left(
\alpha\right)  }\left(  x\right)  dx$ is pseudo-positive. For $k>0,$ by
(\ref{powermmuk2}) and (\ref{eqY}) we obtain
\[
\int r^{-k}d\mu_{k,1}^{\alpha}=2\sqrt{\pi}\int_{0}^{1}r^{k+1}\left(
1-r^{\alpha}\right)  dr=\frac{2\sqrt{\pi}\alpha}{\left(  k+2\right)  \left(
\alpha+k+2\right)  }.
\]
It follows that $w^{\left(  \alpha\right)  }\left(  x\right)  dx$ satisfies
the summability condition (\ref{maincond}), so we can apply our cubature
formula to this kind of measures.

On the other hand, there exist pseudo-positive measures which do not satisfy
the summability condition (\ref{maincond}):

\begin{proposition}
Let $w\left(  re^{it}\right)  :=P\left(  re^{it}\right)  $ for $0\leq r<1$ and
$w\left(  re^{it}\right)  :=0$ for $r\geq1$ where $P\left(  x\right)  $ is
given by (\ref{eqPoisson}). Then $d\mu:=w\left(  x\right)  dx$ is a
pseudo-positive, non-negative moment measure which does not satisfy the
summability condition (\ref{maincond}).
\end{proposition}

\begin{proof}
It follows from (\ref{powermmuk2}) for $k\geq1$
\[
\int r^{-k}d\mu_{k,1}=\int_{0}^{\infty}w_{k,1}\left(  r\right)  \cdot
r^{d-1}dr=2\sqrt{\pi}\int_{0}^{1}r^{k+1}dr=\frac{2\sqrt{\pi}}{\left(
k+2\right)  },
\]
so we see that the summability condition (\ref{maincond}) is not fulfilled.
\end{proof}

Next we compute explicitly the error in Section \ref{Serror} for the function
$w^{\left(  \alpha\right)  }\left(  x\right)  $ with $\alpha=2.$

\begin{theorem}
Let $d\mu:=w^{\left(  2\right)  }\left(  x\right)  dx$, and let $\sigma
^{\left(  s\right)  }$ be the polyharmonic Gau\ss-Jacobi measure of order $s$.
Then for every $f\in C^{2s}\left(  \mathbb{R}^{d}\right)  $ the error
$E_{s}\left(  f\right)  =\int f\left(  x\right)  d\mu-\int f\left(  x\right)
d\sigma^{\left(  s\right)  }$ can be estimated by
\[
\frac{\sqrt{\pi}}{\left(  2s\right)  !}\sum_{k=0}^{\infty}\sup_{0<\xi
<1}\left|  \frac{d^{2s}}{dt^{2s}}\left[  f_{k,l}\left(  \sqrt{t}\right)
t^{-\frac{1}{2}k}\right]  \left(  \xi\right)  \right|  \frac{s!\left(
s+k+1\right)  !}{\left(  2s+k+1\right)  !}\frac{\left(  s+1\right)  !\left(
s+k\right)  !}{\left(  2s+k+2\right)  !}.
\]
\end{theorem}

\begin{proof}
For $k\geq0$ we obtain by (\ref{eqLaplacemu}) and (\ref{eqPoisson}) the
equality
\[
\int_{0}^{1}p\left(  t\right)  d\mu_{k,1}=2\sqrt{\pi}\int_{0}^{1}p\left(
r\right)  \left(  1-r^{2}\right)  r^{2k+1}dr,
\]
and clearly $\mu_{k,2}=0.$ Let $\psi$ be the usual one defined by $\psi\left(
t\right)  =t^{2}.$ Then $\mu_{k,1}^{\psi}$ can be computed by
\begin{equation}
\int_{0}^{1}p\left(  t\right)  d\mu_{k,1}^{\psi}=\int_{0}^{1}p\left(
t^{2}\right)  d\mu_{k,1}=2\sqrt{\pi}\int_{0}^{1}p\left(  r^{2}\right)  \left(
1-r^{2}\right)  r^{2k+1}dr. \label{integrall}%
\end{equation}
According to Theorem \ref{PropError} we have only to compute $\int_{0}^{R^{2}%
}\left|  Q_{k,l}^{s}\left(  t\right)  \right|  ^{2}d\mu_{kl}^{\psi}$ where
$Q_{k,l}^{s}\left(  t\right)  $ is the $s$-th orthogonal polynomial with
leading coefficient $1.$ The substitution $t=r^{2} $ in the integral
(\ref{integrall}) yields $\int_{0}^{1}p\left(  t\right)  d\mu_{k,1}^{\psi
}=\sqrt{\pi}\int_{0}^{1}p\left(  t\right)  \left(  1-t\right)  t^{k}dt.$ The
substitution $t:=\frac{1}{2}\left(  x+1\right)  $ and $dt=\frac{1}{2}dx$ shows
that
\[
\int_{0}^{1}p\left(  t\right)  d\mu_{k,1}^{\psi}=\frac{\sqrt{\pi}}{2^{k+2}%
}\int_{-1}^{1}p\left(  \frac{1}{2}\left(  x+1\right)  \right)  \left(
1-x\right)  \cdot\left(  1+x\right)  ^{k}dx.
\]
Let $P_{s}^{\left(  \alpha,\beta\right)  }\left(  x\right)  $ be the Jacobi
polynomial of degree $s$ (see \cite[p. 30]{gautschi}) normalized with
$P_{s}^{\left(  \alpha,\beta\right)  }\left(  1\right)  =\binom{s+\alpha}{s}.$
They are orthogonal with respect to the measure $dw_{\alpha,\beta}:=\left(
1-x\right)  ^{\alpha}\left(  1+x\right)  ^{\beta}dx$ for $\alpha>-1,\beta>-1.$
It is known that the leading coefficient $k_{s}$ of $P_{s}^{\left(
\alpha,\beta\right)  }\left(  t\right)  $ is equal to $k_{s}=2^{-s}%
\binom{2s+\alpha+\beta}{s}$. Further
\[
h_{s}:=\int_{-1}^{1}\left|  P_{s}^{\left(  \alpha,\beta\right)  }\left(
x\right)  \right|  ^{2}dw_{\alpha,\beta}=\frac{2^{a+\beta+1}}{2s+\alpha
+\beta+1}\frac{\Gamma\left(  s+\alpha+1\right)  \Gamma\left(  s+\beta
+1\right)  }{s!\Gamma\left(  s+\alpha+\beta+1\right)  }.
\]
Define $\widetilde{P}_{s}^{\left(  \alpha,\beta\right)  }\left(  t\right)
:=P_{s}^{\left(  \alpha,\beta\right)  }\left(  2t-1\right)  $ for $0\leq
t\leq1.$ Then $\widetilde{P}_{s}^{\left(  1,k\right)  }$ are orthogonal
polynomials for $\mu_{k,1}^{\psi}$ since
\begin{equation}
\int\widetilde{P}_{s}^{\left(  1,k\right)  }\widetilde{P}_{m}^{\left(
1,k\right)  }d\mu_{k,1}^{\psi}=\frac{\sqrt{\pi}}{2^{k+2}}\int_{-1}^{1}%
P_{s}^{\left(  1,k\right)  }\left(  x\right)  P_{m}^{\left(  1,k\right)
}\left(  x\right)  \left(  1-x\right)  \left(  1+x\right)  ^{k}dx.
\label{letzt}%
\end{equation}
The leading coefficient $\widetilde{k_{s}}$ of $\widetilde{P}_{s}^{\left(
1,k\right)  }\left(  t\right)  =P_{s}^{\left(  1,k\right)  }\left(
2t-1\right)  $ is equal to $2^{s}k_{s}$, so $\widetilde{k_{s}}=\binom
{2s+k+1}{s}.$ Thus $Q_{k,1}^{s}\left(  t\right)  :=\frac{1}{\widetilde{k_{s}}%
}\widetilde{P}_{s}^{\left(  1,k\right)  }$ and we obtain from (\ref{letzt})
\[
\int Q_{k,1}^{s}\left(  t\right)  ^{2}d\mu_{k,1}^{\psi}=\frac{1}%
{\widetilde{k_{s}}^{2}}\frac{\sqrt{\pi}}{2^{k+2}}\int_{-1}^{1}\left|
P_{s}^{\left(  1,k\right)  }\left(  x\right)  \right|  ^{2}\left(  1-x\right)
\left(  1+x\right)  ^{k}dx.
\]
Now this is equal to
\[
\binom{2s+k+1}{s}^{-1}\frac{\sqrt{\pi}}{2^{k+2}}\frac{2^{k+1+1}}{2s+k+2}%
\frac{\Gamma\left(  s+2\right)  \Gamma\left(  s+k+1\right)  }{s!\Gamma\left(
s+k+2\right)  }%
\]
giving
\[
\int Q_{k,1}^{s}\left(  t\right)  ^{2}d\mu_{k,1}^{\psi}=\frac{s!\left(
s+k+1\right)  !}{\left(  2s+k+1\right)  !}\frac{\left(  s+1\right)  !\left(
s+k\right)  !}{\left(  2s+k+2\right)  !}\sim\frac{1}{k^{2s+2}}%
\]
for large $k.$
\end{proof}

\subsection{The summability condition}

We are now turning back to the general situation. The next result shows that
the spectrum of the component measures $\sigma_{k,l}$ is contained in the
spectrum of the representation measure $\mu$.

\begin{theorem}
Let $\sigma_{k,l}$ be non-negative measures on $\left[  0,\infty\right)  $. If
the functional $T:\mathbb{C}\left[  x_{1},x_{2},...,x_{d}\right]
\rightarrow\mathbb{C}$ defined by (\ref{defTTT}) possesses a representing
moment measure $\mu$ with compact support then
\[
\sigma_{k,l}(\left\{  \left|  x\right|  ^{2}\right\}  )\leq\max_{\theta
\in\mathbb{S}^{d-1}}\left|  Y_{k,l}\left(  \theta\right)  \right|
\cdot\left|  x\right|  ^{k}\cdot\left|  \mu\right|  \left(  \left|  x\right|
^{2}\mathbb{S}^{d-1}\right)
\]
for any $x\in\mathbb{R}^{d}$ where $\left|  \mu\right|  $ is the total
variation and $\left|  x\right|  ^{2}\mathbb{S}^{d-1}=\{\left|  x\right|
^{2}\theta:\theta\in\mathbb{S}^{d-1}\}.$
\end{theorem}

\begin{proof}
Let the support of $\mu$ be contained in $B_{R}.$ Let $x_{0}\in\mathbb{R}^{d}
$ be given. For every univariate polynomial $p\left(  t\right)  $ with
$p\left(  \left|  x_{0}\right|  ^{2}\right)  =1$ we have
\begin{align*}
\sigma_{k,l}\left(  \{\left|  x_{0}\right|  ^{2}\}\right)   &  \leq\int
_{0}^{\infty}p\left(  r^{2}\right)  d\sigma_{k,l}\leq\int_{\mathbb{R}^{d}%
}\left|  p(\left|  x\right|  ^{2})Y_{k,l}\left(  x\right)  \right|  d\left|
\mu\right| \\
&  \leq\max_{\theta\in\mathbb{S}^{d-1}}\left|  Y_{k,l}\left(  \theta\right)
\right|  \int_{\mathbb{R}^{d}}\left|  p(\left|  x\right|  ^{2})\right|
\left|  x\right|  ^{k}d\left|  \mu\right|  .
\end{align*}
Now choose a sequence of polynomials $p_{m}$ with $p_{m}\left(  \left|
x_{0}\right|  ^{2}\right)  =1$ which converges on $\left[  0,R\right]  $ to
the function $f$ defined by $f\left(  \left|  x_{0}\right|  ^{2}\right)  =1$
and $f\left(  t\right)  =0$ for $t\neq\left|  x_{0}\right|  ^{2}.$ Since
$\left|  \mu\right|  $ has support in $B_{R}$ Lebesgue's convergence theorem
shows that
\[
\sigma_{k,l}\left(  \{\left|  x_{0}\right|  ^{2}\}\right)  \leq\max_{\theta
\in\mathbb{S}^{d-1}}\left|  Y_{k,l}\left(  \theta\right)  \right|
\int_{\mathbb{R}^{d}}\left|  f\left(  x\right)  \right|  \left|  x\right|
^{k}d\left|  \mu\right|  .
\]
The last implies our statement.
\end{proof}

The following result shows that the summability condition is sometimes
equivalent to the existence of a pseudo-positive representing measure:

\begin{corollary}
Let $d=2.$ Let $\sigma_{k,l}$ be non-negative measures on $\left[
0,\infty\right)  $ and assume that they have disjoint and at most countable
supports. If the functional $T:\mathbb{C}\left[  x_{1},x_{2}\right]
\rightarrow\mathbb{C}$ defined by (\ref{defTTT}) possesses a representing
moment measure with compact support then
\[
\sum_{k=0}^{\infty}\sum_{l=1}^{a_{k}}\int_{0}^{\infty}r^{-k}d\sigma
_{k,l}\left(  r\right)  <\infty\text{.}%
\]
\end{corollary}

\begin{proof}
Let $\Sigma_{k,l}$ be the support set of $\sigma_{k,l}.$ The last theorem
shows that $\sigma_{k,l}\left(  \left\{  0\right\}  \right)  =0,$ hence
$0\notin\Sigma_{k,l}.$ Moreover it tells us that
\[
\int_{0}^{\infty}r^{-k}d\sigma_{k,l}\left(  r\right)  \leq\max_{\theta
\in\mathbb{S}^{d-1}}\left|  Y_{k,l}\left(  \theta\right)  \right|  \cdot
\sum_{r\in\Sigma_{k,l}}\left|  \mu\right|  \left(  r\mathbb{S}^{d-1}\right)
.
\]
Since $d=2$ we know that $\max_{\theta\in\mathbb{S}^{d-1}}\left|
Y_{k,l}\left(  \theta\right)  \right|  \leq1.$ Hence
\[
\sum_{k=0}^{\infty}\sum_{l=1}^{a_{k}}\int_{0}^{\infty}r^{-k}d\sigma
_{k,l}\left(  r\right)  \leq\sum_{k=0}^{\infty}\sum_{l=1}^{a_{k}}\sum
_{r\in\Sigma_{k,l}}\left|  \mu\right|  \left(  r\mathbb{S}^{d-1}\right)
\leq\left|  \mu\right|  \left(  \mathbb{R}^{d}\right)
\]
where the last inequality follows from the fact that $\Sigma_{k,l}$ are
pairwise disjoint.
\end{proof}

Recall that the converse of the last theorem holds under the additional
assumption that the supports of all $\sigma_{k,l}$ are contained in some
interval $\left[  0,R\right]  .$

\begin{theorem}
\label{Tnonexistence}There exists a functional $\ T:\mathbb{C}\left[
x_{1},x_{2},...,x_{d}\right]  \rightarrow\mathbb{C}$ which is pseudo-positive
definite but does not possess a pseudo-positive representing measure.
\end{theorem}

\begin{proof}
Let $\sigma$ be a non-negative measure over $\left[  0,R\right]  .$ Let
$f\in\mathbb{C}\left[  x_{1},..,x_{d}\right]  $ and let $f_{k,l}$ be the
Laplace-Fourier coefficients of $f.$ By Proposition \ref{PropTTT} it is clear
that
\[
T\left(  f\right)  :=\int_{0}^{R}f_{1,1}\left(  r\right)  r^{-1}d\sigma\left(
r\right)
\]
is pseudo-positive definite. We take now for $\sigma$ the Dirac functional at
$r=0$. Suppose that $T$ has a signed representing measure $\mu$ which is
pseudo-positive. Then the component measure $\mu_{11}$ is non-negative, and it
is defined by the equation
\[
\int_{0}^{\infty}h\left(  r\right)  d\mu_{11}\left(  r\right)  :=\int
_{\mathbb{R}^{n}}h\left(  \left|  x\right|  \right)  Y_{11}\left(  x\right)
d\mu
\]
for any continuous function $h:\left[  0,\infty\right)  \rightarrow\mathbb{C}$
with compact support. Take now $h\left(  r\right)  =r^{2}.$ Then by
Proposition \ref{pseudopos}
\[
\int_{0}^{\infty}r^{2}d\mu_{11}\left(  r\right)  =\int_{\mathbb{R}^{n}}\left|
x\right|  ^{2}Y_{11}\left(  x\right)  d\mu=T\left(  \left|  x\right|
^{2}Y_{11}\left(  x\right)  \right)  =0.
\]
It follows that $\mu_{11}$ has support $\left\{  0\right\}  .$ On the other
hand, if we take a sequence of functions $h_{m}\in C_{c}\left(  \left[
0,\infty\right)  \right)  $ such that $h_{m}\rightarrow1_{\left\{  0\right\}
},$ then we obtain
\[
\mu_{11}\left(  \left\{  0\right\}  \right)  =\lim_{m\rightarrow\infty}%
\int_{\mathbb{R}^{n}}h_{m}\left(  \left|  x\right|  \right)  Y_{11}\left(
x\right)  d\mu.
\]
But $h_{m}\left(  \left|  x\right|  \right)  Y_{11}\left(  x\right)  $
converges to the zero-function, and Lebesgue's theorem shows that $\mu
_{11}\left(  \left\{  0\right\}  \right)  =0,$ so $\mu_{11}=0.$ This is a
contradiction since
\[
\int_{0}^{\infty}1d\mu_{11}\left(  r\right)  =\int_{\mathbb{R}^{n}}%
Y_{11}\left(  x\right)  d\mu=T\left(  Y_{11}\right)  =\int_{0}^{R}%
1d\sigma\left(  r\right)  =1.
\]
The proof is complete.
\end{proof}

\section{Concluding Remarks \label{Sfinal}}

One important feature of the polyharmonic Gau\ss--Jacobi cubature which
deserves to be discussed is its numerical significance. At first glance, one
may object that the cubature needs the knowledge of the Laplace-Fourier
coefficients of the function $f:\mathbb{R}^{n}\rightarrow\mathbb{C},$ and
these are based on integrals as well. However, if one works with polynomials,
the Gau\ss\ decomposition (\ref{Gaussdecomposition}) can be constructed by an
efficient differentiation algorithm, see \cite{AxRa95}. Decomposing the
harmonic polynomials $h_{j}$ according to our fixed orthonormal basis
$Y_{k,l}\left(  x\right)  $ one obtains from (\ref{Gaussdecomposition}) the
expansion
\[
f\left(  x\right)  =\sum_{k=0}^{\deg f}\sum_{l=1}^{a_{k}}p_{k,l}(\left|
x\right|  ^{2})Y_{k,l}\left(  x\right)  ,
\]
where $p_{k,l}$ are uniquely determined univariate polynomials. For the
Gau\ss-Jacobi quadrature there exists suitable software to compute the weights
$\alpha_{1}^{k,l},...,\alpha_{s}^{k,l}$ and the nodes $r_{1}^{\left(
k,l\right)  },...,r_{s}^{\left(  k,l\right)  }$ of the measures $\sigma
_{k,l}^{\left(  s\right)  }$. Then
\[
T^{\left(  s\right)  }\left(  f\right)  :=\sum_{k=0}^{\deg f}\sum_{l=1}%
^{a_{k}}\sum_{j=1}^{s}\alpha_{j}^{k,l}p_{k,l}\left(  (r_{j}^{\left(
k,l\right)  })^{2}\right)  .
\]
In practice one also has to bound the number of spherical harmonics
$Y_{k,l}\left(  x\right)  $ in the formula. Our main results Theorem
\ref{ThmMain} and Theorem \ref{TStieltjes} show that by increasing the number
$k$ and the number $s$ the algorithm remains stable.

What concerns the polyharmonic Gau\ss--Jacobi cubature in the annulus and the
strip, we note that recently efficient algorithms for finding Gau\ss
--Jacobi--Markov quadratures for Chebyshev systems have been studied in
\cite{MRW96}.

A second point to be made clear, is the motivation why we choose the space of
polyharmonic functions of order $s$ as the exactness space for the
multivariate generalization. One main reason is the fact that recently
polyharmonic functions have shown to be an efficient tool in approximation
theory and more generally, in mathematical analysis, see e.g. \cite{ACL83},
\cite{kounchev92}, \cite{kounchev98}, \cite{Koun00}, and \cite{Ligo88}.
Another motivation stems from potential theory: two non-negative measures
$\mu$ and $\nu$ with compact supports are \emph{gravitationally equivalent}
if
\begin{equation}
\int_{\mathbb{R}^{n}}h\left(  x\right)  d\mu\left(  x\right)  =\int
_{\mathbb{R}^{n}}h\left(  x\right)  d\nu\left(  x\right)  \label{grav}%
\end{equation}
for all harmonic functions $h$ defined on a neighborhood of the supports of
the measures. If $\mu$ and $\nu$ are gravitationally equivalent\emph{\ }then
they produce the same potential outside of their support. Graviequivalent
measures are very important in inverse problems in Geophysics and Geodesy, and
a new mathematical area has grown extensively during the last two decades or
so under the title Quadrature Domains, see the comprehensive survey and
references in \cite{gustafssonshapiro}, as well as \cite{zidarov}.

In analogy to (\ref{grav}) one could define two (generally speaking, signed)
measures $\mu$ and $\nu$ as polyharmonically equivalent of order $s$ if
(\ref{grav}) holds for all polyharmonic functions of order $s$ in a
neighborhood of their support. Similar notions of equivalence have been
developed by L. Ehrenpreis \cite{ehrenpreis} in the form of a generalized
balayage. To our knowledge the equivalence of two measures (and more
generally, distributions) with respect to the solutions of an elliptic
operator has been for the first time rigorously formulated and studied in the
case of non--negativity in \cite{schulzewildenhain}.

Let us remark that the polyharmonic Gau\ss--Jacobi measure $\nu^{\left(
s\right)  }$ which we have introduced in the present paper is related to the
concept of ''mother body'' (a non--negative measure $\nu$ satisfying
(\ref{grav}) for a given $\mu$, and having minimal support) in the theory of
Quadrature Domains, cf. \cite{gustafssonshapiro}, and Remark
\ref{Rminimalsupport}, as well as \cite{kounchev87} and \cite{kounchev93}.

ACKNOWLEDGMENT. Both authors acknowledge the support of the Institutes
Partnership project with the Alexander von Humboldt Foundation.

Author's addresses:

1. Ognyan Kounchev, Institute of Mathematics and Informatics, Bulgarian
Academy of Sciences, 8 Acad. G. Bonchev Str., 1113 Sofia, Bulgaria;

e--mail: kounchev@math.bas.bg, kounchev@math.uni--duisburg.de

2. Hermann Render, Departamento de Matem\'{a}ticas y Computati\'{o}n,
Universidad de la Rioja, Edificio Vives, Luis de Ulloa, s/n. 26004
Logro\~{n}o, Spain; e-mail: render@math.uni-duisburg.de; herender@dmc.unirioja.es
\end{document}